\documentclass[preprint,12pt,authoryear]{elsarticle}

\usepackage{a4wide}
\usepackage{amssymb}
\usepackage{amsfonts}

\usepackage{amsmath,amsthm,amsfonts,amssymb,color}

\usepackage{multicol}

\usepackage{natbib}

\usepackage{multirow} 
\usepackage{float}
\usepackage{graphicx}
\usepackage{longtable}  

\usepackage{amssymb}
\usepackage{amsthm}

\usepackage[authoryear]{natbib}

\begin{document}

\begin{frontmatter}

\title{Lie symmetry analysis of the nonlinear generalized heat equation for varying cross-section geometry} 


\author{Targyn A. Nauryz$^{1,2,3}$}

\address{(1)International School of Economics, Kazakh-British Technical University, Almaty, Kazakhstan\\
(2)Institute of Mathematics and Mathematical Modeling, Almaty, Kazakhstan\\
(3)School of Digital Technologies, Narxoz University, Almaty, Kazakhstan\\
Email: targyn.nauryz@gmail.com}

\begin{abstract}
We study the nonlinear generalized heat equation $C(u)u_t=\frac{1}{z^{\nu}}\left(K(u)z^{\nu}u_z\right)_z$,
where $C(u)$ and $K(u)$ are temperature-dependent thermal coefficients and $\nu>0$ is a geometric parameter describing the varying cross-section geometry. By applying the classical Lie symmetry method, we derive the determining equations and perform a complete classification of the admitted Lie point symmetries according to the functional dependence between $C(u)$ and $K(u)$. The analysis shows that the symmetry structure splits naturally into two principal cases: $C(u)/K(u)$ non-constant and $C(u)/K(u)=\beta$ constant. In the first case, only the basic symmetries are admitted for arbitrary coefficients, whereas additional generators appear under special compatibility relations. In the second case, the equation can be transformed to a linear radial heat equation by the substitution $v=\int K(u)du$, yielding an extended symmetry algebra. For each case, we construct the infinitesimal generators, commutator tables, one-parameter transformation groups, and corresponding invariant reductions. Invariant and similarity solutions are obtained and then specialized to several physically relevant subclasses, including power-law, exponential-type, and linear constitutive coefficients. The results provide a unified symmetry-based model for the analysis of generalized nonlinear heat equations in non-Cartesian geometries.
\end{abstract}

\begin{keyword}
    generalized heat equation, Lie symmetry method, infinitesimals, commutator, invariant solution, similarity solution.
\end{keyword}
\end{frontmatter}

\section{Introduction}

Nonlinear heat equations with variable thermal coefficients arise naturally in many areas of applied mathematics, heat transfer, phase transition theory, and engineering applications. In particular, when the heat capacity and thermal conductivity depend on the temperature, the governing equations acquire a nonlinear structure that makes the derivation of exact solutions substantially more difficult than in the constant-coefficient setting. At the same time, such models provide a more realistic description of thermal processes in heterogeneous and temperature-sensitive media.

The Lie group method is one of the most effective analytical tools for the investigation of nonlinear differential equations. Its theoretical foundations, including infinitesimal transformations, prolongations, invariant surface conditions, and symmetry reductions, are well established in the classical monographs \cite{1}-\cite{3}, \cite{15}, \cite{17}, \cite{18}. Further algebraic background for symmetry structures and Lie algebras may be found in \cite{4}. These works provide the general framework for determining admitted symmetry groups, constructing similarity variables, and reducing partial differential equations to ordinary differential equations.

Lie symmetry analysis has been widely applied to heat-type and nonlinear evolution equations. For example, in \cite{14}  the nonlinear heat equation $u_t=u_{xx}+f(u)$ is studied and symmetry reductions together with exact solutions for several classes of nonlinearities $f(u)$ are obtained. The classical group classification of the nonlinear heat equation may be traced back to \cite{19}, while a recent revisitation of the general nonlinear heat equation $u_t=(K(u)u_x)_x$ was presented in \cite{20}. Closely related diffusion-convection models were analyzed in \cite{23}, who studied the equation $u_t=(D(u)u_x)_x-K^{\prime}(u)u_x$, and derived group classifications and symmetry reductions according to the form of the transport coefficients. Variable-coefficient nonlinear diffusion-convection equations of the form $f(x)u_t=(g(x)D(u)u_x)_x-k^{\prime}(u)u_x$ were investigated in \cite{32}, while symmetry-based exact solutions for Burgers-type equations, $u_t+\alpha u u_x+\beta u_{xx}=0$ and their generalizations were considered in \cite{16}, \cite{21}, \cite{22}. These studies demonstrate that Lie symmetry methods are particularly well suited for identifying special coefficient structures that enhance the admitted symmetry algebra and allow exact reductions.

Symmetry methods have also been successfully extended to generalized and fractional models. In \cite{5}, the heat equation is considered through modified one-parameter local point transformations. Therefore, in \cite{7} Lie symmetries of a generalized modified heat equation is investigated. Lie and related transformation methods are applied to several fractional models, including coupled fractional KdV systems \cite{9}, time-fractional Boussinesq equations \cite{10}, Boussinesq-type systems \cite{11}, and other diffusion, wave, and telegraph equations \cite{8}. These contributions show that symmetry analysis remains a powerful and flexible approach not only for classical parabolic equations, but also for more general nonlinear and fractional evolution equations.

Another important direction concerns nonlinear boundary-value and free-boundary problems. In \cite{12}, nonlinear boundary value problems are studied with Lie symmetries, while a complete Lie symmetry classification for a class of $(1+2)$-dimensional reaction-diffusion-convection equations carried out in \cite{13}. These works are especially relevant because, in many applications, the usefulness of symmetry methods depends not only on the governing equation itself, but also on the compatibility of the admitted symmetries with boundary and interface conditions.

In recent years, generalized heat equations and Stefan-type models with nonlinear thermal coefficients have attracted increasing attention in connection with phase-change processes and heat transfer in non-Cartesian geometries. In \cite{24}, the nonlinear heat-diffusion equation $C(u)u_t=(K(u)u_x)_x$ is studied, performed a Lie symmetry classification according to the functional relation between $C(u)$ and $K(u)$, and derived invariant solutions for several constitutive laws. This equation may be viewed as a one-dimensional Cartesian counterpart of the model considered in the present paper.

The study of generalized heat equations is also strongly motivated by problems arising in electrical-contact theory and Stefan-type thermal processes. A two-phase Stefan problem for a generalized heat equation is inverstigated with special functions methods \cite{25}. In \cite{26}, two-phase Stefan problem for a generalized heat equation with nonlinear thermal coefficients is studied. The mathematical modeling of the heat process in closure electrical contacts with a heat source is considered \cite{27}. In \cite{28}, the author analyzed the thermal process of arc erosion with current-carrying heating effect in a temperature gradient and in [29], an exact solution to a Stefan-type problem for a generalized heat equation with Thomson effect derived. A nonlinear Stefan problem with internal Joule heat source in cylindrical geometry ($\nu=1$) by fixed point theory is studied \cite{30}. Finally, in \cite{31}, the mathematical modeling of a heat process in a cylindrical domain with nonlinear thermal coefficients and a heat source on the axis considered \cite{31}. These works show that generalized heat equations with nonlinear thermal parameters and variable geometry are not only mathematically interesting, but also of direct relevance for the modeling of electrical-contact phenomena, melting, vaporization, and related moving-boundary processes.

Motivated by the above developments, in this paper we study the nonlinear generalized heat equation
\begin{equation}\label{eq1}
    C(u)u_t=\frac{1}{z^{\nu}}\left(K(u)z^{\nu}u_z\right)_z,
\end{equation}
where $\nu>0$ is a geometric parameter characterizing the varying cross-section geometry. This model extends the Cartesian nonlinear heat-diffusion equation by incorporating an explicit radial-type geometric term. In particular, the term $\frac{1}{z^{\nu}}\left(K(u)z^{\nu}u_z\right)_z$ distinguishes the present equation from the standard one-dimensional case and makes the symmetry structure essentially dependent on the geometry.

The main objective of this work is to perform a Lie point symmetry analysis of the above generalized heat equation, to classify the admitted infinitesimal generators according to the functional relation between $C(u)$ and $K(u)$, and to construct the corresponding invariant and similarity solutions. A key feature of the analysis is that the symmetry classification splits naturally into two principal cases depending on whether the ratio $C(u)/K(u)$ is constant or non-constant.

In addition to the abstract classification, we also consider several constitutive subclasses of physical interest, including power-law, exponential-type, and linear conductivity/capacity pairs. These examples make it possible to write the generators explicitly and to derive concrete invariant solutions, thereby linking the general Lie symmetry framework with nonlinear thermal models of practical interest.

The paper is organized as follows. In Section 2, we derive the determining equations and classify the admitted Lie point symmetries of the generalized heat equation. In Section 3, we construct the corresponding commutator tables and one-parameter transformation groups. Section 4 is devoted to invariant reductions and the derivation of similarity solutions associated with the admitted generators. In Section 5, we specialize the general results to several physically relevant classes of nonlinear coefficients. Finally, the last section summarizes the main results and outlines possible directions for future work.

\section{Lie group symmetries of the governing equation}
We consider the heat conduction equation \eqref{eq1} with temperature depended coefficients $C(u)$ and $K(u)$, which equivalent to
\begin{equation}\label{eq2}
    C(u)u_t=K^{\prime}(u)u_z^2+K(u)u_{zz}+\frac{\nu}{z}K(u)u_z,
\end{equation}
where $\nu>0$ is a geometric parameter characterizing the underlying spatial symmetry. 

Our aim is to determine the continuous transformation groups leaving this equation invariant and to classify the admitted Lie point symmetries depending on the functional forms of the coefficients $C(u)$ and $K(u)$.

To this end, we consider a one-parameter local Lie group of point transformations acting on the variables $(z,t,u)$ with the group parameter $\lambda$
\begin{equation}\label{eq3}
    \begin{cases}
        \Bar{z}=Z(z,t,u,\lambda)\\
        \Bar{t}=T(z,t,u,\lambda)\\
        \Bar{u}=U(z,t,u,\lambda),
    \end{cases}
\end{equation}
The associated infinitesimal transformations of \eqref{eq3} are given by
\begin{equation}\label{eq4}
    \begin{cases}
        \Bar{z}=z+\lambda\xi(z,t,u)+o(\lambda^2)\\
        \Bar{t}=t+\lambda\tau(z,t,u)+o(\lambda^2)\\
        \Bar{u}=u+\lambda\eta(z,t,u)+o(\lambda^2),
    \end{cases}
\end{equation}
where 
$$\xi(z,t,u)=\frac{\partial Z}{\partial \lambda}\bigg|_{\lambda=0},\quad \tau(z,t,u)=\frac{\partial T}{\partial \lambda}\bigg|_{\lambda=0},\quad \eta(z,t,u)=\frac{\partial U}{\partial \lambda}\bigg|_{\lambda=0}.$$

If we define 
\begin{equation}\label{eq5}
    G(z,u,u_z,u_t,u_{zz})= C(u)u_t-K^{\prime}(u)u_z^2-K(u)u_{zz}-\frac{\nu}{z}K(u)u_z,
\end{equation}
then the invariance of equation \eqref{eq1} under the transformation group \eqref{eq4} is equivalent to the infinitesimal invariance condition
\begin{equation}\label{eq6}
    Y^{(2)}(G)\big|_{G=0}=0,
\end{equation}
where $Y^{(2)}$ is the second prolongation of the infinitesimal generator
\begin{equation}\label{eq7}
    Y=\xi(z,t,u)\partial_z+\tau(z,t,u)\partial_t+\eta(z,t,u)\partial_u,
\end{equation}
which is stated by
\begin{equation}\label{eq8}
    Y^{(2)}=Y+\eta_z^{(1)}\partial_{u_z}+\eta_{t}^{(1)}\partial_{u_t}+\eta_{zz}^{(2)}\partial_{u_{zz}}+\eta_{zt}^{(2)}\partial_{u_{zt}}+\eta_{tz}^{(2)}\partial_{u_{tz}}+\eta_{tt}^{(2)}\partial_{u_{tt}}. 
\end{equation}
If $G_t=G_{u_{zt}}=G_{u_{tz}}=G_{u_{tt}}=0$, we can derive 
\begin{equation}\label{eq9}
    Y^{(2)}(G)=\xi G_z+\eta G_u+\eta_z^{(1)}G_{u_z}+\eta_{t}^{(1)}G_{u_t}+\eta_{zz}^{(2)}G_{u_{zz}},
\end{equation}
where
\begin{equation}\label{eq10}
    \eta_z^{(1)}=\eta_z+(\eta_u-\xi_z)u_z-\tau_z u_t-\xi_u u_z^2-\tau_u u_z u_t,
\end{equation}
\begin{equation}\label{eq11}
    \eta_t^{(1)}=\eta_t+(\eta_u-\tau_t)u_t-\xi_t u_z-\tau_{u}u_t^2-\xi_u u_z u_t,
\end{equation}
\begin{equation}\label{eq12}
    \begin{split}
        \eta_{zz}^{(2)}=&\eta_{zz}+(2\eta_{zu}-\xi_{zz})u_z-\tau_{zz}u_t+(\eta_{uu}-2\xi_{zu})u_z^2-2\tau_{zu}u_z u_t\\
        &+(\eta_u-2\xi_z)u_{zz}-2\tau_z u_{zt}-\xi_{uu}u_z^3-\tau_{uu}u_z^2 u_t-3\xi_u u_z u_{zz}\\
        &-\tau_u u_t u_{zz}-2\tau_u u_z u_{zt}.
    \end{split}
\end{equation}

Applying the condition \eqref{eq6} to equation \eqref{eq2} we state 

\begin{equation}\label{eq13}
    \begin{split}
        &\eta\left(C^{\prime}(u)u_t-K^{\prime\prime}(u)u_z^2-K^{\prime}(u)u_{zz}-\frac{\nu}{z}K^{\prime}(u)u_z\right)-\eta^{(1)}_z\left(2K^{\prime}(u)u_z+\frac{\nu}{z}K(u)\right)\\
        &+C(u)\eta^{(1)}_t-K(u)\eta_{zz}^{(2)}+\xi\frac{\nu}{z^2}K(u)u_z=0.
    \end{split}
\end{equation}

From the full invariance condition \eqref{eq2}, terms involving $u_{zt},\;u_z u_{zt},\;u_t u_{zz},\;u_z u_{zz},\;u_{t}^2,\;u_z u_t$ must be eliminated. Hence, their coefficients force

\begin{equation}\label{eq14}
    \xi_u=0,\quad\quad \tau_z=0,\quad\quad \tau_u=0.
\end{equation}

which lead to

\begin{equation}\label{eq15}
    \xi=\xi(z,t),\quad\quad \tau=\tau(t).
\end{equation}

Collecting terms without any derivatives of $u$ we have

\begin{equation}\label{eq16}
    C(u)\eta_t-K(u)\eta_{zz}-\frac{\nu}{z}K(u)\eta_z=0.
\end{equation}

Thus, from the coefficients of $u_z,\;u_{z}^2,\;u_{zz}$ in equation \eqref{eq13} we obtain 

\begin{equation}\label{eq17}
    2K^{\prime}(u)\eta_z+C(u)\xi_t+K(u)(2\eta_{zu}-\xi_{zz})+\nu K(u)\left(\frac{\xi_z}{z}-\frac{\xi}{z^2}\right)=0,
\end{equation}

\begin{equation}\label{eq18}
    \eta\left(\frac{C^{\prime}(u)K^{\prime}(u)}{C(u)}-K^{\prime\prime}(u)\right)+K^{\prime}(u)(2\xi_z-\eta_u-\tau^{\prime})-K(u)\eta_{uu}=0,
\end{equation}

\begin{equation}\label{eq19}
    \eta\left(\frac{C^{\prime}(u)K(u)}{C(u)}-K^{\prime}(u)\right)+K(u)(2\xi_z-\tau^{\prime})=0.
\end{equation}

Solving the equation \eqref{eq19} for $\eta$, we have

\begin{equation}\label{eq20}
    \eta=F(u)\left(\tau^{\prime}-2\xi_{z}\right),\quad\text{where}\quad F(u)=\left(\frac{C^{\prime}(u)}{C(u)}-\frac{K^{\prime}(u)}{K(u)}\right)^{-1},
\end{equation}

and plugging \eqref{eq20} into equation \eqref{eq16}, we get

\begin{equation}\label{eq21}
    \frac{C(u)}{K(u)}\left(2\xi_{zt}-\tau^{\prime\prime}\right)=2\xi_{zzz}+\frac{2\nu}{z}\xi_{zz}.
\end{equation}

Because the ratio $C(u)/K(u)$ is a function of $u$ alone, while the remaining terms in equation \eqref{eq21} involve only the independent variables $z$ and $t$, equation \eqref{eq21} can be satisfied only in certain distinct cases. This leads naturally to the consideration of two possibilities: either $C(u)/K(u)$ is non-constant, or $C(u)/K(u)=\beta$ constant.

\section{$C(u)/K(u)$ non constant}
If we consider the case when $C(u)/K(u)$ is depended on $u$, then from \eqref{eq21} we provide the equations
\begin{equation}\label{eq22}
    2\xi_{zt}-\tau^{\prime\prime}=0,
\end{equation}
\begin{equation}\label{eq23}
    \xi_{zzz}+\frac{\nu}{z}\xi_{zz}=0,
\end{equation}
applying the substitution $w(z,t)=\xi_{zz}(z,t)$ to the equation \eqref{eq23} which leads to  
\begin{itemize}
    \item $\nu\neq 1,2:\quad \xi_{zz}=a(t)z^{-\nu}$ for general case and integrating two times we have $\xi(z,t)=p(t)z^{2-\nu}+b(t)z+\mu(t),\quad p(t)=\frac{a(t)}{(1-\nu)(2-\nu)}$. From \eqref{eq22} we derive  
    \begin{equation}\label{eq24}
        2(2-\nu)p^{\prime}(t)z^{1-\nu}+2b^{\prime}(t)-\tau^{\prime\prime}=0.
    \end{equation}
    Therefore the coefficients of the independent functions of $z$ must vanish:
    $$p^{\prime}(t)=0,\quad\quad 2b^{\prime}(t)-\tau^{\prime\prime}=0.$$
    So $\rho(t)\equiv\rho$ and $b(t)=\frac{\tau^{\prime}}{2}+\alpha$, thus
    \begin{equation}\label{eq25}
        \xi(z,t)=\rho z^{2-\nu}+\left(\frac{\tau^{\prime}}{2}+\alpha\right)z+\gamma(t),
    \end{equation}
    where $\rho,\alpha$ are constants. Hence, from \eqref{eq20} we can obtain
    \begin{equation}\label{eq26}
        \eta(z,t,u)=-2F(u)((2-\nu)\rho z^{1-\nu}+\alpha).
    \end{equation}
    According to the equation \eqref{eq17}, we have
    \begin{equation}\label{eq27}
        \begin{split}
            &(1-\nu)\rho z^{-\nu}\left(-4(2-\nu)\frac{K^{\prime}(u)}{C(u)}F(u)-4(2-\nu)\frac{K(u)}{C(u)}F^{\prime}(u)-2\frac{K(u)}{C(u)}\right)\\
            &+\frac{\tau^{\prime\prime}(t)}{2}z+\gamma^{\prime}(t)-\nu\frac{K(u)}{C(u)}\frac{\gamma(t)}{z^2}=0.
        \end{split}
    \end{equation}
    Since $\frac{C(u)}{K(u)}$ is non constant, equivalently $\frac{K(u)}{C(u)}$ is nonconstant, the terms containing $\frac{K(u)}{C(u)}$ and $\frac{K^{\prime}(u)}{C(u)}$ must separate from the purely $z,t$ dependent terms. Therefore
    \begin{equation}\label{eq28}
        (1-\nu)\rho z^{-\nu}\left(-4(2-\nu)\frac{K^{\prime}(u)}{C(u)}F(u)-4(2-\nu)\frac{K(u)}{C(u)}F^{\prime}(u)-2\frac{K(u)}{C(u)}\right)=0,
    \end{equation}
    \begin{equation}\label{eq29}
        \tau^{\prime\prime}(t)=0,\quad\quad \gamma(t)=0.
    \end{equation}
    From \eqref{eq28}, we can consider
    \begin{equation}\label{eq30}
        \rho=0 \quad \text{or}\quad 2(2-\nu)\frac{K^{\prime}(u)}{C(u)}F(u)+2(2-\nu)\frac{K(u)}{C(u)}F^{\prime}(u)+\frac{K(u)}{C(u)}=0.
    \end{equation}
    The equation \eqref{eq29} implies 
    \begin{equation}\label{eq31}
        \tau(t)=\sigma t+\delta,
    \end{equation}
    with $\sigma,\delta$ are constants. Thus, the generic-case infinitesimals are
    \begin{equation}\label{eq32}
        \begin{cases}
            \xi(z)=\rho z^{2-\nu}+\left(\frac{\sigma}{2}+\alpha\right)z,\\
            \tau(t)=\sigma t+\delta,\\
            \eta(z,u)=-2F(u)((2-\nu)\rho z^{1-\nu}+\alpha).
        \end{cases}
    \end{equation}

    \item $\nu=1:\;\xi_{zz}=a(t)z^{-1}$ for cylindrical case, analogously solving \eqref{eq22},\eqref{eq23} and \eqref{eq29}, we have
          \begin{equation}\label{eq33}
              \begin{cases}
                  \xi(z)=\rho(z\ln z-z)+\left(\frac{\sigma}{2}+\alpha\right)z,\\
                  \tau(t)=\sigma t+\delta,\\
                  \eta(z,u)=-2F(u)(\rho 
                  \ln z+\alpha),
              \end{cases}
          \end{equation}

    \item $\nu=2:\;\xi_{zz}=a(t)z^{-2}$  for spherical case, with similar approach, we can obtain
        \begin{equation}\label{eq34}
            \begin{cases}
                \xi(z)=-\rho \ln z+\left(\frac{\sigma}{2}+\alpha\right)z,\\
                \tau(t)=\sigma t+\delta,\\
                \eta(z,u)=-2F(u)\left(\frac{\rho}{z}+\alpha\right).
            \end{cases}
        \end{equation}
\end{itemize}
with four parameters $\rho,\alpha,\sigma,\delta$.

The validity of the obtained infinitesimals is completed by enforcing equation \eqref{eq18}, which is equivalent to the relation below
\begin{equation}\label{eq35}
    \left((2-\nu)\rho z^{1-\nu}+\alpha\right)\left[-F(u)\left(\frac{C^{\prime}(u)}{C(u)}K^{\prime}(u)-K^{\prime\prime}(u)\right)+K^{\prime}(u)+K^{\prime}(u)F^{\prime}(u)+K(u)F^{\prime\prime}(u)\right]=0.
\end{equation}

It should be noted that the parameters $\sigma$ and $\delta$ exist for arbitrary functions $K(u)$ and $C(u)$, while relation \eqref{eq35} shows that the existence of the parameters $\alpha$ and $\rho$ is determined by the specific choice of $K(u)$ and $C(u)$. Hence, three distinct cases must be considered:

\begin{itemize}
    \item[a)] $\rho=0$ and $\alpha=0$. In this case, \eqref{eq32} gives $\xi(z)=\frac{\sigma}{2}z$, $\tau(t)=\sigma t+\delta$ and $\eta(z,u)=0$. This means that the dependent variable $u$ remains unchanged, while the independent variables $z$ and $t$ undergo scaling-type transformations. In particular, the parameter $\sigma$ generates a simultaneous dilation in space and time, whereas $\delta$ corresponds to a translation in time. Hence, the admitted infinitesimal generators are
    \begin{equation}\label{eq36}
        Y_1=\frac{z}{2}\partial_z+t\partial_t,\quad\quad Y_2=\partial_t.
    \end{equation}
    It should be emphasized that, in contrast to the Cartesian case \cite{24}, the operator $\partial_z$ is generally not admitted here. Indeed, the governing equation contains the explicit term $\nu/z$, which introduces direct dependence on the spatial variable $z$. Therefore, the equation is not invariant under translations in $z$ for $\nu>0$. As a consequence, a constant spatial shift cannot be included among the admitted symmetries.

    \item[b)] $\rho=0$ and $\alpha\neq 0$. Now we consider the case $\xi(z)=\left(\frac{\sigma}{2}+\alpha\right)z$, $\tau(t)=\sigma t+\delta$ and $\eta(z,u)=-2\alpha F(u)$,
    where $F(u)$ is defined by \eqref{eq20}.

    Since $\rho=0$, condition \eqref{eq30} is satisfied identically. Therefore, it remains to impose condition \eqref{eq35}, which yields
    $$-F(u)\left(\frac{C^{\prime}(u)}{C(u)}K^{\prime}(u)-K^{\prime\prime}(u)\right)+K^{\prime}(u)+K^{\prime}(u)F^{\prime}(u)+K(u)F^{\prime\prime}(u)=0,$$
    diving this equation by $K(u)$ and using identity
    $$\left(\frac{K^{\prime}(u)}{K(u)}\right)^{\prime}=\frac{K^{\prime\prime}(u)}{K(u)}-\left(\frac{K^{\prime}(u)}{K(u)}\right)^2,$$
    together with the definition of $F(u)$, the previous relation can be rewritten in the form
    $$F(u)\left(\frac{K^{\prime}(u)}{K(u)}\right)^{\prime}+\left(\frac{K^{\prime}(u)}{K(u)}\right)F^{\prime}(u)+F^{\prime\prime}(u)=0,$$
    which equivalent to
    $$\left(F(u)\frac{K^{\prime}(u)}{K(u)}\right)^{\prime}=-F^{\prime\prime}(u).$$
    Integrating with respect to $u$, we arrive at
    \begin{equation}\label{eq37}
        (F(u)K(u))^{\prime}=A K(u),
    \end{equation}
    where $A$ is an arbitrary constant. Hence,
    \begin{equation}\label{eq38}
        F(u)=\frac{A\int K(u)du+B}{K(u)},
    \end{equation}
    with $B$ being another integration constant. Substituting this expression into the defining relation \eqref{eq20} we obtain the compatibility condition between the coefficient functions $C(u)$ and $K(u)$. At this stage, two subcases must be considered.

    If $A\neq 0$, from \eqref{eq38}, it follows that 
    \begin{equation}\label{eq39}
        C(u)=D K(u)\left(A\int K(u)du+B\right)^{1/A},
    \end{equation}
    where $D$ is a nonzero constant.

    In the case $A=0$, we can obtain
    \begin{equation}\label{eq40}
        C(u)=D K(u)\exp\left(\frac{1}{B}\int K(u)du\right).
    \end{equation}
    Thus, in both subcases, the existence of an additional symmetry requires a special functional relation between $C(u)$ and $K(u)$. Under these conditions, the extra infinitesimal generator takes the form
    \begin{equation}\label{eq41}
        Y_3=z\partial_z-2\frac{A\int K(u)du+B}{K(u)}\partial_u.
    \end{equation}
    Therefore, besides the basic generators admitted for arbitrary $C(u)$ and $K(u)$, the equation admits the additional operator $Y_3$ whenever the coefficients satisfy either relation \eqref{eq39} or relation \eqref{eq40}.

    \item[c)] $\rho\neq 0$ and $\alpha=0$. We have $\xi(z)=\rho z^{2-\nu}+\frac{\sigma}{2}z$, $\tau(t)=\sigma t+\delta$ and $\eta(z,u)=-2(2-\nu)\rho F(u)z^{1-\nu}$. Then multiplying both sides of \eqref{eq30} by $C(u)$, we get
    $$2(2-\nu)(K(u)F(u))^{\prime}+K(u)=0.$$
    Hence
    $$(K(u)F(u))^{\prime}=-\frac{K(u)}{2(2-\nu)}$$
    and integrating, we obtain
    $$K(u)F(u)=-\frac{1}{2(2-\nu)}\int K(u)du +M,$$
    where $M$ is constant. Therefore
    \begin{equation}\label{eq42}
        F(u)=\frac{2(2-\nu)M-\int K(u)du}{2(2-\nu)K(u)}.
    \end{equation}
    Using \eqref{eq20}, for $\nu\neq 2$ and taking into account that $F$ is given by \eqref{eq20}, it then follows that
    \begin{equation}\label{eq43}
        C(u)=\frac{NK(u)}{\left(2(2-\nu)M-\int K(u)du\right)^{\frac{2}{2-\nu}}}
    \end{equation}
    with $N$ constant. 
    It should be noted that the obtained expression for \(C(u)\) represents a particular case of the general form of \(C(u)\) given in \eqref{eq40}, corresponding to the choice $M=B$, $N=(2(2-\nu))^{2(2-\nu)}D$ and $A=-\frac{1}{2(2-\nu)}$. Therefore, equation \eqref{eq35} is satisfied. Then the new generator corresponding to $\rho$ is
    \begin{equation}\label{eq44}
        Y_4=z^{2-\nu}\partial_z-\frac{2(2-\nu)M-\int K(u)du}{K(u)}z^{1-\nu}\partial_u.
    \end{equation}
    For cylindrical case $\nu=1$, the generator \eqref{eq44} takes the form
    \begin{equation}\label{eq45}
        Y_4=(z\ln z-z)\partial_z-\ln z\frac{2M-\int K(u)du}{K(u)}\partial_u,
    \end{equation}
    and for spherical case $\nu=2$, equation \eqref{eq44} becomes
    \begin{equation}\label{eq46}
        Y_4=-\ln z\partial_z
    \end{equation}
\end{itemize}
Therefore, the case $\rho\neq 0$ and $\alpha\neq 0$ does not produce a new independent structure, but rather combines the generators obtained in cases (b) and (c). The full generator is given by the linear combination of $Y_1$, $Y_2$, $Y_3$ and $Y_4$ provided the corresponding coefficient restrictions are satisfied simultaneously.

In summary, the infinitesimal generators for the general case $\nu\neq 1,2$ take the following form

\begin{equation}\label{eq47}
    \begin{cases}
        Y_1=\frac{z}{2}\partial_z+t\partial_t,\\
        Y_2=\partial_t,\\
        Y_3=z\partial_z-2\frac{A\int K(u)du+B}{K(u)}\partial_u,\\
        Y_4=z^{2-\nu}\partial_z-\frac{2(2-\nu)M-\int K(u)du}{K(u)}z^{1-\nu}\partial_u,
    \end{cases}
\end{equation}
where $F(u)$ is given by \eqref{eq20}.

The Lie bracket is defined in the usual way by
$$\left[Y_i, Y_j\right]=Y_i Y_j - Y_j Y_i$$
with the generators interpreted as differential operators. The generators $Y_1$ and $Y_2$ are admitted for arbitrary functions $C(u)$ and $K(u)$, whereas the operators $Y_3$ and $Y_4$ appear only under additional restrictions on the coefficient functions. The four-dimensional Lie algebra generated by the infinitesimal operators is characterized by the commutation relations listed in Table \ref{tab1}.
\begin{table}[ht]
\centering
\caption{Commutator table for non constant case of $\frac{C(u)}{K(u)}$}
\scriptsize
\renewcommand{\arraystretch}{2} 
\begin{tabular}{|p{1cm}|p{1.5cm}|p{1.5cm}|p{1.5cm}|p{1.5cm}|}
\hline $\left[\;,\;\right]$ & $Y_1$ & $Y_2$ & $Y_3$ & $Y_4$\\
\hline $Y_1$ & 0 & $-Y_2$ & 0 & $\frac{1-\nu}{2}Y_4$ \\
\hline $Y_2$ & $Y_2$ & 0 & 0 & 0 \\
\hline $Y_3$ & 0 & 0 & 0 & $(1-\nu)Y_4$\\
\hline $Y_4$ & $-\frac{1-\nu}{2}Y_4$ & 0 & $-(1-\nu)Y_4$ & 0\\
\hline
\end{tabular}\label{tab1}
\end{table}
Then, for each infinitesimal generator, we obtain the corresponding one-parameter Lie group of transformations by solving the first-order system
$$\frac{d\Bar{z}}{d\lambda}=\xi(\Bar{z},\Bar{t},\Bar{u},\lambda),\quad\quad\frac{d\Bar{t}}{d\lambda}=\tau(\Bar{z},\Bar{t},\Bar{u},\lambda),\quad\quad \frac{d\Bar{u}}{d\lambda}=\eta(\Bar{z},\Bar{t},\Bar{u},\lambda),$$
subject to the initial conditions $\Bar{z}=z$, $\Bar{t}=t$, $\Bar{u}=u$ at $\lambda=0$.

The resulting one-parameter groups associated with each infinitesimal generator are listed below:

$$G_1:=\{\Bar{z}=ze^{\lambda/2},\quad \Bar{t}=te^{\lambda},\quad \Bar{u}=u\},$$
$$G_2:=\{\Bar{z}=z,\quad \Bar{t}=t+\lambda,\quad \Bar{u}=u\},$$
\begin{equation*}
    G_3=\begin{cases}
        \Bar{z}=ze^{\lambda},\quad \Bar{t}=t,\quad \Bar{u}=J^{-1}\left(\left(J(u)+\frac{B}{A}\right)e^{-2A\lambda}-\frac{B}{A}\right),& A\neq 0\\
        \Bar{z}=ze^{\lambda},\quad \Bar{t}=t,\quad \Bar{u}=J^{-1}\left(J(u)-2B\lambda\right), & A=0,
    \end{cases}
\end{equation*}
$$G_4=\left\{\Bar{z}=(z^{1-\nu}+(\nu-1)\lambda)^{\frac{1}{\nu-1}},\quad \Bar{t}=t,\quad \Bar{u}=J^{-1}\left(2(2-\nu)M+(J(u)-2(2-\nu)M)e^{z^{1-\nu}\lambda}\right)\right\}$$
and hence, for $\nu=1$ and $\nu=2$ we can construct the following Lie groups
$$G_4^{(\nu=1)}:=\left\{\Bar{z}=\exp(1+(\ln z-1)e^{\lambda}),\quad \Bar{t}=t,\quad \Bar{u}=J^{-1}\left(2M+(J(u)-2M)e^{-\lambda \ln z}\right)\right\},$$
$$G_4^{(\nu=2)}:=\bigg\{\Bar{z}=li^{-1}(li(z)-\lambda),\quad \Bar{t}=t,\quad \Bar{u}=u\bigg\},$$
where $J(u)=\int K(u)du$ is assumed to be invertible and $li(z)=\int_0^z\frac{dz}{\ln z}$.

\section{$C(u)/K(u)=\beta$ constant}

If we consider that the ratio $\frac{C(u)}{K(u)}$ is a constant such that $C(u)=\beta K(u)$, then equations \eqref{eq16}, \eqref{eq17}, \eqref{eq18}, \eqref{eq19} and \eqref{eq20} are equivalent to the equations
\begin{equation}\label{eq48}
    \beta\eta_t-\eta_{zz}-\frac{\nu}{z}\eta_z=0,
\end{equation}
\begin{equation}\label{eq49}
    2\frac{K^{\prime}(u)}{K(u)}\eta_z+\beta\xi_t+2\eta_{zu}-\xi_{zz}+\nu\left(\frac{\xi_z}{z}-\frac{\xi}{z}\right)=0,
\end{equation}
\begin{equation}\label{eq50}
    \eta\left(\frac{(K^{\prime}(u))^2}{K^2(u)}-\frac{K^{\prime\prime}(u)}{K(u)}\right)+\frac{K^{\prime}(u)}{K(u)}(2\xi_z-\eta_u-\tau^{\prime})-\eta_{uu}=0,
\end{equation}
\begin{equation}\label{eq51}
    2\xi_z-\tau^{\prime}=0.
\end{equation}
Integrating the equation \eqref{eq51} respect to $z$, one can provide
\begin{equation}\label{eq52}
    \xi(z,t)=\frac{1}{2}\tau^{\prime}z+a(t),
\end{equation}
with $a(t)$ is an arbitrary function. Thus, from \eqref{eq50} we can state
\begin{equation}\label{eq53}
    \frac{\partial}{\partial u}\left(\eta_u+\eta\frac{K^{\prime}(u)}{K(u)}\right)=0,
\end{equation}
which equivalent to 
\begin{equation}\label{eq54}
    \eta=-\left(\eta_u+g(z,t)\right)\frac{K(u)}{K^{\prime}(u)},
\end{equation}
and solving \eqref{eq54} for $\eta$, it follows that 
\begin{equation}\label{eq55}
    \eta(z,t,u)=-g(z,t)\frac{\int K(u)du}{K(u)}.
\end{equation}
Substituting \eqref{eq52} and \eqref{eq54} into equation \eqref{eq49} and solving respect to $g(z,t)$, we have
\begin{equation}\label{eq56}
    g(z,t)=\frac{\beta}{8}\tau^{\prime\prime}z^2+\frac{\beta}{2}a^{\prime}(t)z+\frac{\nu a(t)}{2z}+d(t).
\end{equation}
From \eqref{eq48} and \eqref{eq54} one can obtain $\beta g_t-g_{zz}-\frac{\nu}{z}g_z=0$ and plugging \eqref{eq56} into this equation, it takes the following form

\begin{equation}\label{eq57}
    \frac{\beta^2}{8}\tau^{\prime\prime\prime}z^2+\frac{\beta^2}{2}a^{\prime\prime}(t)z+\beta d^{\prime}(t)-\frac{\beta(1+\nu)}{4}\tau^{\prime\prime}+\nu\left(\frac{\nu}{2}-1\right)\frac{a(t)}{z^3}=0.
\end{equation}
Therefore, we consider 
\begin{equation}\label{eq58}
    \tau^{\prime\prime\prime}(t)=0,\quad\quad a^{\prime\prime}(t)=0,\quad\quad d^{\prime}(t)-\frac{\beta(1+\nu)}{4}\tau^{\prime\prime}=0,\quad\quad \nu\left(\frac{\nu}{2}-1\right)a(t)=0.
\end{equation}
A crucial restriction comes from the last condition, $\nu\left(\frac{\nu}{2}-1\right)a(t)=0$.
Because $\nu>0$, the behavior of $a(t)$ depends on whether $\nu\neq 2$ or $\nu=2$. Indeed, if $\nu=2$, then the factor $\nu\left(\frac{\nu}{2}-1\right)$ does not vanish, and therefore one must have $a(t)=0$. In contrast, when $\nu=2$, this coefficient becomes zero identically, and the function $a(t)$ is no longer forced to vanish. In that case, it may survive as a nontrivial linear function of $t$.

Therefore, the further symmetry classification naturally splits into two separate cases.
\begin{itemize}
    \item[1)] If $\nu\neq 2$, then $a(t)=0$ and the second, third conditions in \eqref{eq58} become
    \begin{equation}\label{eq59}
        \tau(t)=\mu t^2+\gamma t+m,\quad\quad d(t)=\frac{1+\nu}{2}\mu t+k,
    \end{equation}
    from \eqref{eq56} we have
    \begin{equation}\label{eq60}
        g(z,t)=\frac{\beta}{4}\mu z^2+\frac{1+\nu}{2}\mu t+k.
    \end{equation}
    In summary, applying \eqref{eq52}, \eqref{eq55} and \eqref{eq59}, we get the following infinitesimals
    \begin{equation}\label{eq61}
        \begin{cases}
            \xi(z,t)=\mu tz+\frac{\gamma}{2}z\\
            \tau(t)=\mu t^2+\gamma t+m\\
            \eta(z,t,u)=-\frac{\int K(u)du}{K(u)}\left(\frac{\beta}{4}\mu z^2+\frac{1+\nu}{2}\mu t+k\right),
        \end{cases}
    \end{equation}
    then the generic radial case $\nu\neq 2$ admits the four-parameter Lie algebra generated by
    \begin{equation}\label{eq62}
        \begin{cases}
            \widetilde{Y}_1=tz\partial_z + t^2\partial_t -\left(\frac{\beta}{4}z^2+\frac{1+\nu}{2}t\right)\frac{\int K(u)du}{K(u)}\partial_u,\\
            \widetilde{Y}_2=\frac{z}{2}\partial_z+t\partial_t,\\
            \widetilde{Y}_3=\partial_t,\\
            \widetilde{Y}_4=-\frac{\int K(u)du}{K(u)}\partial_u,
        \end{cases}
    \end{equation}
    with constants $\mu,\gamma,m,k,\beta$. 
    
    Hence, the commutator table of Lie algebra constructed from infinitesimal generators \eqref{eq62} is presented in Table \ref{tab2}. 
    \begin{table}[ht]
        \centering
        \caption{Commutator table for constant case $C(u)=\beta K(u)$ and $\nu\neq 2$}
        \scriptsize
        \renewcommand{\arraystretch}{2} 
        \begin{tabular}{|p{1cm}|p{1.8cm}|p{1.5cm}|p{1.9cm}|p{1cm}|}
        \hline $\left[\;,\;\right]$ & $\widetilde{Y}_1$ & $\widetilde{Y}_2$ & $\widetilde{Y}_3$ & $\widetilde{Y}_4$\\
        \hline $\widetilde{Y}_1$ & 0 & $-\widetilde{Y}_1$ & $-2\widetilde{Y}_2-\frac{1+\nu}{2}\widetilde{Y}_4$ & 0 \\
        \hline $\widetilde{Y}_2$ & $\widetilde{Y}_1$ & 0 & $-\widetilde{Y}_3$ & 0 \\
        \hline $\widetilde{Y}_3$ & $2\widetilde{Y}_2+\frac{1+\nu}{2}\widetilde{Y}_4$ & $\widetilde{Y}_3$ & 0 & 0\\
        \hline $\widetilde{Y}_4$ & 0 & 0 & 0 & 0\\
        \hline
        \end{tabular}\label{tab2}
    \end{table}
    The Lie groups of transformations corresponding to each infinitesimal generators \eqref{eq62} can be defined as follows:
    $$L_1=\left\{\Bar{z}=\frac{z}{1-\lambda t},\;\;\Bar{t}=\frac{t}{1-\lambda t},\;\;\Bar{u}=J^{-1}\left(J(u)(1-\lambda t)^{\frac{1+\nu}{2}}\exp\left(-\frac{\beta\lambda z^2}{4(1-\lambda t)}\right)\right)\right\},$$
    $$L_2=\left\{\Bar{z}=ze^{\lambda/2},\;\;\Bar{t}=te^{\lambda t},\;\;\Bar{u}=u\right\},$$
    $$L_3=\left\{\Bar{z}=z,\;\; \Bar{t}=t+\lambda,\;\;\Bar{u}=u\right\},$$
    $$L_4=\left\{\Bar{z}=z,\;\;\Bar{t}=t,\;\;\Bar{u}=J^{-1}\left(J(u)e^{-\lambda}\right)\right\},$$
    with $J(u)=\int K(u)du$ which considered as invertible function.
    \item[2)] Hence, analogously for $\nu=2$, we can list the following infinitesimals
    \begin{equation}\label{eq63}
        \begin{cases}
            \xi(z,t)=\mu zt+\frac{\gamma}{2}z+ht+\omega,\\
            \tau(t)=\mu t^2+\gamma t+m\\
            \eta(z,t,u)=-\frac{\int K(u)du}{K(u)}\left(\frac{\beta}{4}\mu z^2+\frac{\beta}{2}hz+\frac{ht+\omega}{z}+\frac{3}{2}\mu t+k\right),
        \end{cases}
    \end{equation}
    thus the six infinitesimal generators are
    \begin{equation}\label{eq64}
        \begin{cases}
            \widehat{Y}_1=tz\partial_z+t^2\partial_t-\left(\frac{\beta}{4}z^2+\frac{3}{2}t\right)\frac{\int K(u)du}{K(u)}\partial_u,\\
            \widehat{Y}_2=\frac{z}{2}\partial_z+t\partial_t,\\
            \widehat{Y}_3=\partial_t,\\
            \widehat{Y}_4=t\partial_z-\left(\frac{\beta}{2}z+\frac{t}{z}\right)\frac{\int K(u)du}{K(u)}\partial_u,\\
            \widehat{Y}_5=\partial_z-\frac{1}{z}\frac{\int K(u)du}{K(u)}\partial_u,\\
            \widehat{Y}_6=-\frac{\int K(u)du}{K(u)}\partial_u,
        \end{cases}
    \end{equation}
    with constants $\beta,\mu,\gamma,h,\omega,m$ and $k$.
    These operators form a six-dimensional Lie algebra whose structure is described by the commutation relations, see Table \ref{tab3}. The commutator table is useful for identifying the algebraic properties of the admitted symmetry group, in particular, for distinguishing central elements and for constructing optimal systems of subalgebras in subsequent reduction procedures.
    \begin{table}[ht]
        \centering
        \caption{Commutator table for constant case $C(u)=\beta K(u)$ and $\nu=2$}
        \scriptsize
        \renewcommand{\arraystretch}{2} 
        \begin{tabular}{|p{1cm}|p{1.5cm}|p{1.5cm}|p{1.5cm}|p{1.5cm}|p{1.5cm}|p{1cm}|}
        \hline $\left[\;,\;\right]$ & $\widehat{Y}_1$ & $\widehat{Y}_2$ & $\widehat{Y}_3$ & $\widehat{Y}_4$ & $\widehat{Y}_5$ & $\widehat{Y}_6$\\
        \hline $\widehat{Y}_1$ & 0 & $-\widehat{Y}_1$ & $-2\widehat{Y}_2-\frac{3}{2}\widehat{Y}_6$ & 0 & $-\widehat{Y}_4$ & 0 \\
        \hline $\widehat{Y}_2$ & $\widehat{Y}_1$ & 0 & $-\widehat{Y}_3$ & $-\frac{1}{4}\widehat{Y}_4$ & $-\frac{1}{2}\widehat{Y}_5$ & 0 \\
        \hline $\widehat{Y}_3$ & $2\widehat{Y}_2+\frac{3}{2}\widehat{Y}_6$ & $\widehat{Y}_3$ & 0 & $\widehat{Y}_5$ & 0 & 0\\
        \hline $\widehat{Y}_4$ & 0 & $-\frac{1}{2}\widehat{Y}_4$ & $-\widehat{Y}_5$ & 0 & $-\frac{\beta}{2}\widehat{Y}_6$ & 0\\
        \hline $\widehat{Y}_5$ & $\widehat{Y}_4$ & $\frac{1}{2}\widehat{Y}_5$ & 0 & $\frac{\beta}{2}\widehat{Y}_6$ & 0 & 0\\
        \hline $\widehat{Y}_6$ & 0 & 0 & 0 & 0 & 0 & 0\\
        \hline
        \end{tabular}\label{tab3}
    \end{table}
    From this table one can see that $\Tilde{Y}_6$ commutes with all basis elements, and therefore it spans the center of the Lie algebra. The generators $\Tilde{Y}_4$ and $\Tilde{Y}_5$ form a nontrivial substructure coupled through $\Tilde{Y}_6$, while $\Tilde{Y}_1,\;\Tilde{Y}_2,\;\Tilde{Y}_3$ describe the projective-dilation-translation part of the symmetry algebra in the variables $z$ and $t$.

    Next, integrating the Lie equations associated with each infinitesimal generator in \eqref{eq64}, we obtain the corresponding one-parameter Lie groups of transformations.
    $$\Tilde{L}_1=\left\{\Bar{z}=\frac{z}{1-\lambda t},\;\;\Bar{t}=\frac{t}{1-\lambda t},\;\; \Bar{u}=J^{-1}\left(J(u)(1-\lambda t)^{3/2}\exp\left(-\frac{\beta\lambda z^2}{4(1-\lambda t)}\right)\right)\right\},$$
    $$\Tilde{L}_2=\left\{\Bar{z}=ze^{\lambda/2},\;\;\Bar{t}=te^{\lambda},\;\;\Bar{u}=u\right\},$$
    $$\Tilde{L}_3=\left\{\Bar{z}=z,\;\;\Bar{t}=t+\lambda,\;\;\Bar{u}=u\right\},$$
    $$\Tilde{L}_4=\left\{\Bar{z}=z+\lambda t,\;\;\Bar{t}=t,\;\;\Bar{u}=J^{-1}\left(J(u)\frac{z}{z+\lambda t}\exp\left(-\frac{\beta}{2}\lambda z-\frac{\beta}{4}\lambda^2 t\right)\right)\right\},$$
    $$\Tilde{L}_5=\left\{\Bar{z}=z+\lambda,\;\;\Bar{t}=t,\;\;\Bar{u}=J^{-1}\left(J(u)\frac{z}{z+\lambda}\right)\right\},$$
    $$\Tilde{L}_6=\left\{\Bar{z}=z,\;\;\Bar{t}=t,\;\;\Bar{u}=J^{-1}\left(J(u)e^{-\lambda}\right)\right\},$$
    here we assume that $J(u)=\int K(u)du$ is an invertible function.
\end{itemize}
 
These finite transformations play a central role in the construction of invariant solutions, since they allow one to map known solutions into new ones and to derive similarity variables for symmetry reduction.

Therefore, the case $\nu=2$ admits a significantly richer symmetry group than the generic radial case. This enlarged symmetry structure makes it possible to construct broader classes of invariant and similarity solutions and provides a stronger algebraic approach for the reduction of the original nonlinear generalized heat equation to ordinary differential equations.

\section{Invariant solutions corresponding to the admitted infinitesimal generators}

Having determined the admitted Lie point symmetries of equation \eqref{eq2}, we now proceed to construct invariant solutions corresponding to each infinitesimal generator. Such solutions are of particular interest because they reduce the number of independent variables and transform the original nonlinear partial differential equation into an ordinary differential equation.

A solution $u=T(x,t)$ is said to be invariant of \eqref{eq2} if the surface $u=T(x,t)$ remains invariant under the one-parameter transformation group generated by $Y$. This requirement is equivalent to the invariance condition
\begin{equation}\label{eq65}
    Y(u-T(x,t))=0\quad\quad\text{on}\quad\quad u=T(x,t).
\end{equation}
This condition leads to a first-order linear equation for $T(z,t)$, or, equivalently, to the characteristic system defining the invariants of the generator. Using these invariants, we derive the corresponding similarity forms of the solution and substitute them into PDE \eqref{eq2} to obtain reduced ordinary differential equations.

Below, we consider the admitted infinitesimal generators one by one and construct the corresponding invariant solutions for each classified case.

\subsection{Infinitesimal generators in non constant case of $\frac{C(u)}{K(u)}$}

\textbf{Cenerator $Y_1$.}  We first analyze the invariant solutions associated with the generator $Y_1$ is given by \eqref{eq47}. Let $u=T(x,t)$ be an invariant solution of equation \eqref{eq2}. Then, by the invariance criterion, $T(x,t)$ must satisfy
\begin{equation}\label{eq66}
    \frac{z}{2}\frac{\partial T}{\partial z}+t\frac{\partial T}{\partial t}=0,
\end{equation}
The corresponding characteristic system is therefore
\begin{equation}\label{eq67}
    \frac{dz}{z/2}=\frac{dt}{t}.
\end{equation}
It follows that $u$ is constant along the characteristic curves, while integration of the first two relations gives the similarity variable
\begin{equation}\label{eq68}
    \eta=\frac{z}{\sqrt{t}}.
\end{equation}
Hence, the invariant solution associated with $Y_1$ can be represented in the form
\begin{equation}\label{eq69}
    u(\eta,t)=\varphi_1(\eta)
\end{equation}
with $\eta$ given by \eqref{eq68}. Substituting this similarity ansatz into equation \eqref{eq2}, we reduce the original partial differential equation to an ordinary differential equation for the profile function $\varphi_1$. After straightforward calculation, we obtain
\begin{equation}\label{eq70}
    K(\varphi_1)\varphi_1^{\prime\prime}+K^{\prime}(\varphi_1)\left(\varphi^{\prime}_1\right)^2+\frac{\nu}{\eta}K(\varphi_1)\varphi_1^{\prime}+\frac{\eta}{2} C(\varphi_1)\varphi_1^{\prime}=0,
\end{equation}
where the prime denotes differentiation with respect to $\eta$. This equation may also be rewritten in form as
\begin{equation}\label{eq71}
    \left(K(\varphi_1)\varphi_1^{\prime}\right)^{\prime}+\left(\frac{\nu}{\eta}+\frac{\eta}{2}\frac{C(\varphi_1)}{K(\varphi_1)}\right)K(\varphi_1)\varphi_1^{\prime}=0.
\end{equation}
Introducing the auxiliary function $f(\eta)=K(\varphi_1)\varphi_1^{\prime}$, we arrive at the first-order equation
$$f^{\prime}(\eta)+\left(\frac{\nu}{\eta}+\frac{\eta}{2}\frac{C(\varphi_1)}{K(\varphi_1)}\right)f(\eta)=0.$$
Its formal solution is
\begin{equation}\label{eq72}
    f(\eta)=C_1\eta^{-\nu}\exp\left(-\frac{1}{2}\int\eta\frac{C(\varphi_1)}{K(\varphi_1)}d\eta\right),
\end{equation}
where $C_1$ is an arbitrary constant. Consequently, the invariant profile $\varphi_1$ satisfies the integral relation
\begin{equation}\label{eq73}
    \varphi_1(\eta)=C_2+C_1\int\frac{\eta^{-\nu}}{K(\varphi_1(\eta))}\exp\left(-\frac{1}{2}\int \eta\frac{C(\varphi_1(\eta))}{K(\varphi_1(\eta))}d\eta\right)d\eta.
\end{equation}
with $C_2$ being another integration constant. 

\textbf{Generator $Y_2$.} We next consider the generator $Y_2$ defined by \eqref{eq42}. If $u=T(x,t)$ is an invariant solution of equation \eqref{eq2}, then the invariance condition implies
\begin{equation}\label{eq74}
    \frac{\partial T}{\partial t}=0.
\end{equation}
Hence, any $Y_2$-invariant solution is stationary and depends only on the spatial variable $z$. Therefore, we set
\begin{equation}\label{eq75}
    u(z,t)=\varphi_2(z).
\end{equation}
Substituting this expression into equation \eqref{eq2}, we obtain
\begin{equation*}
    \left(K(\varphi_2)\varphi_2^{\prime}\right)^{\prime}+\frac{\nu}{\eta}K(\varphi_2)\varphi_2^{\prime}=0,
\end{equation*}
or equivalently,
\begin{equation}\label{eq76}
    \frac{1}{z^{\nu}}\frac{d}{dz}\left(z^{\nu}K(\varphi_2)\varphi_2^{\prime}\right)=0.
\end{equation}
Integrating once gives
$$z^{\nu}K(\varphi_2)\varphi_2^{\prime}=C_1,$$
where $C_1$ is a constant. Since $H^{\prime}(\varphi_2)=K(\varphi_2)$, then the above equation takes the form $\frac{d}{dz}H(\varphi_2)=C_1 z^{-\nu}$. A further integration yields
\begin{equation}\label{eq77}
    H(\varphi_2)=\begin{cases}
        \frac{C_1}{1-\nu}z^{1-\nu}+C_2,& \nu\neq 1,\\
        C_1\ln z+C_2,&\nu=1,
    \end{cases}
\end{equation}
where $C_2$ is another integration constant. Consequently, the invariant solution corresponding to $Y_2$ is given by
\begin{equation}\label{eq78}
    u(z,t)=H^{-1}\left( \frac{C_1}{1-\nu}z^{1-\nu}+C_2\right),\quad\quad \nu\neq 1,
\end{equation}
and for the cylindrical case $\nu=1$ we have
\begin{equation}\label{eq79}
    u(z,t)=H^{-1}\left(C_1\ln z+C_2\right).
\end{equation}
Thus, the symmetry generated by $Y_2$ gives rise to the class of time-independent solutions of equation \eqref{eq2}.

\textbf{Generator $Y_3$.} We now turn to the additional generator which can be rewritten in the form
$$Y_3=z\partial_z-2\frac{AJ(u)+B}{K(u)}\partial_u,$$
where $J(u)=\int K(u)du$.
Let $u=T(x,t)$ be invariant under $Y_3$. Then the invariance condition becomes
\begin{equation}\label{eq80}
    z\frac{\partial T}{\partial z}+2\frac{AJ(T)+B}{K(T)}=0.
\end{equation}
It is more convenient to rewrite this relation in terms of the function $J(u)$. Since $\frac{\partial}{\partial z}J(T)=K(T)T_z$, then the above equation is equivalent to
\begin{equation}\label{eq81}
    \frac{dz}{z}=-\frac{dJ(T)}{2(AJ(T)+B)}.
\end{equation}
Hence, $t$ is an invariant, while the second relation yields the similarity representation.

If $A\neq 0$, integration of \eqref{eq81} gives
\begin{equation}\label{eq82}
    z\varphi_3(t)=\left(AJ(u)+B\right)^{-\frac{1}{2A}},
\end{equation}
which is equivalent to
\begin{equation}\label{eq83}
    AJ(u)+B=z^{-2A}\varphi_3(t)^{-2A}.
\end{equation}
Differentiating both sides of \eqref{eq83} with respect to $z$, we have
\begin{equation}\label{eq84}
    K(u)u_z=-2z^{-2A-1}\varphi_3(t)^{-2A},
\end{equation}
and one more differentiation with respect to $z$ implies
\begin{equation}\label{eq85}
    K^{\prime}(u)(u_z)^2+K(u)u_{zz}=2(2A+1)z^{-2A-2}\varphi_3(t)^{-2A}.
\end{equation}
If we differentiate \eqref{eq83} with respect to $t$, we can state
\begin{equation}\label{eq86}
    K(u)u_t=-2Az^{-2A}\varphi_3^{-2A-1}(t)\varphi_3^{\prime}(t).
\end{equation}
Substituting \eqref{eq39} and \eqref{eq82}-\eqref{eq86} into the PDE \eqref{eq2}, we can derive the following equation in terms of $\varphi_3$
\begin{equation}\label{eq87}
    Dz^{-2}\varphi_3^{-2}(t)\cdot\left(-2z^{-2A}\varphi_3(t)^{-2A-1}\varphi_3^{\prime}(t)\right)=2(2A+1-\nu)z^{-2A-2}\varphi_3(t)^{-2A}.
\end{equation}
Solving the equation \eqref{eq87} for $\varphi_3$, we can express
\begin{equation}\label{eq88}
    \varphi_3^2(t)=\frac{D}{QD+2(2A+1-\nu)t}.
\end{equation}
where $Q$ is an arbitrary constant. Finally, from \eqref{eq83} we have
\begin{equation}\label{eq89}
    AJ(u)+B=z^{-2A}\left(Q+\frac{2(2A+1-\nu)}{D}t\right)^A.
\end{equation}
Assuming that $J(u)$ is an invertible, thus the invariant solution is given
\begin{equation}\label{eq90}
    u(z,t)=J^{-1}\left(\frac{z^{-2A}\left(Q+\frac{2(2A+1-\nu)}{D}t\right)^A}{A}-\frac{B}{A}\right).
\end{equation}
In the degenerate case $A=0$, the generator takes the form
\begin{equation}\label{eq91}
    Y_3=z\partial_z-\frac{2B}{K(u)}\partial_u,
\end{equation}
then equation \eqref{eq81} becomes
\begin{equation}\label{eq92}
    \frac{dz}{z}=-\frac{dJ(T)}{2B},
\end{equation}
and equation \eqref{eq82} takes the form
\begin{equation}\label{eq93}
    z\varphi_3(t)=\exp\left(-\frac{J(u)}{2B}\right).
\end{equation}
Differentiating \eqref{eq93} with respect to $z$, we obtain
\begin{equation}\label{eq94}
    K(u)u_z=-2B\varphi_3(t)\exp\left(-\frac{J(u)}{2B}\right),
\end{equation}
and differentiating \eqref{eq93} with respect to $t$, we have
\begin{equation}\label{eq95}
    K(u)u_t=-2Bz\varphi_3^{\prime}(t)\exp\left(-\frac{J(u)}{2B}\right).
\end{equation}
Therefore
\begin{equation}\label{eq96}
    K^{\prime}(u)u_z^2+K(u)u_{zz}=\frac{2B(1-\nu)}{z^2}.
\end{equation}
Substituting \eqref{eq94}, \eqref{eq95} and \eqref{eq96} into the equation \eqref{eq2}, we can derive the following equation
\begin{equation}\label{eq97}
    Dz^{-2}\varphi_3^{-2}(t)\left(-2Bz\varphi_3^{\prime}(t)\exp\left(-\frac{J(u)}{2B}\right)\right)=\frac{2B(1-\nu)}{z^2}.
\end{equation}
Using $e^{-\frac{J(u)}{2B}}=(z\varphi_3(t))^{-1}$, then \eqref{eq97} reduces to
\begin{equation}\label{eq98}
    -D\varphi_3(t)^{-3}\varphi_3^{\prime}(t)=1-\nu.
\end{equation}
Solving equation \eqref{eq98} for $\varphi_3$, one can state the function
\begin{equation}\label{eq99}
    \varphi_3^2(t)=\frac{D}{QD+2(1-\nu)t}.
\end{equation}
Hence the invariant solution is represented by
\begin{equation}\label{eq100}
    u(z,t)=J^{-1}\left(-2B\ln z+B\ln\left(Q+\frac{2(1-\nu)}{D}t\right)\right).
\end{equation}
Therefore, the symmetry generated by $Y_3$ produces a family of logarithmic or power-type similarity solutions, depending on whether $A=0$ or $A\neq 0$.

\textbf{Generator $Y_4$.} Finally, we consider the generator which can be presented in the form
\begin{equation}\label{eq101}
    Y_4=z^{2-\nu}\partial_z-2(2-\nu)F(u)\partial_z,
\end{equation}
where
$$F(u)=\frac{2(2-\nu)M-J(u)}{2(2-\nu)K(u)},\quad\quad J(u)=\int K(u)du,$$
If $u=T(x,t)$ is invariant under $Y_4$, then the invariance condition gives
\begin{equation}\label{eq102}
    -2(2-\nu)F(u)z^{1-\nu}-z^{2-\nu}T_z=0.
\end{equation}
The characteristic system takes the form
\begin{equation}\label{eq103}
    \frac{dz}{z}=\frac{dJ(u)}{2(2-\nu)N-J(u)}.
\end{equation}
Hence integration gives
\begin{equation}\label{eq104}
    2(2-\nu)M-J(u)=z\varphi_4(t).
\end{equation}
Applying $K^{\prime}(u)(u_z)^2+K(u)u_{zz}=0$ and
$$K(u)u_z=-\varphi_4(t),\quad K(u)u_t=-z\varphi_4^{\prime}(t),\quad \frac{C(u)}{K(u)}=\frac{N}{(2(2-\nu)M-J(u))^{\frac{2}{2-\nu}}},$$
the equation \eqref{eq2} can be rewritten
\begin{equation}\label{eq105}
    \frac{Nz\varphi_4^{\prime}(t)}{(z\varphi_4(t))^{\frac{2}{2-\nu}}}=\frac{\nu}{z}\varphi_4(t).
\end{equation}
For $\nu\neq 1,2$, the powers of $z$ are incompatible unless $\varphi_4^{\prime}(t)=0$ and $\varphi_4(t)=0$. Therefore this branch yields only the trivial constant solution
\begin{equation}\label{eq106}
    u(z,t)=J^{-1}(2(2-\nu)M),\quad\quad u(z,t)=\text{const}.
\end{equation}

\subsection{Infinitesimal generators in constant case $C(u)=\beta K(u)$}

We now consider the second principal case arising from the symmetry classification, namely when the ratio of the coefficient functions is constant. Indeed, let us introduce the new dependent variable
\begin{equation}\label{eq107}
    v=J(u)=\int K(u)du,
\end{equation}
where the function $J(u)$ is assumed to be invertible. Then, substituting into equation \eqref{eq2}, and using the relation $C(u)=\beta K(u)$, we obtain
\begin{equation}\label{eq108}
    \beta v_t=v_{zz}+\frac{\nu}{z}v_z.
\end{equation}
This reduction is important for two reasons. First, it shows that the nonlinear model becomes linear after an appropriate change of dependent variable. Second, it allows us to use the symmetry structure of the linear radial heat equation in order to construct invariant and similarity solutions in a systematic way.

Therefore, in this case, our strategy is the following. We first determine the invariant solutions of the linear equation for the transformed variable $v(z,t)$. After that, returning to the original variable is straightforward through the inverse transformation
\begin{equation}\label{eq109}
    u(z,t)=J^{-1}(v(z,t)).
\end{equation}
Consequently, every invariant solution obtained for $v$ immediately generates a corresponding invariant solution of the original nonlinear equation.

\textbf{Generator $\widetilde{Y}_1$.} Let $v=v(z,t)$ be a solution invariant under the action of $\widetilde{Y}_1$ given by \eqref{eq62}. Then, by the invariance criterion, it must satisfy
$$\widetilde{Y}_1(v-v(z,t))\big|_{v=v(z,t)}=0,$$
which leads to the first-order linear equation
\begin{equation}\label{eq110}
    tzv_z+t^2v_t+\left(\frac{\beta}{4}z^2+\frac{1+\nu}{2}t\right)v=0.
\end{equation}
To solve this equation, we consider its characteristic system
\begin{equation}\label{eq111}
    \frac{dz}{tz}=\frac{dt}{t^2}=-\frac{dv}{\left(\frac{\beta}{4}z^2+\frac{1+\nu}{2}t\right)v}.
\end{equation}
From the first two relations, we obtain
\begin{equation}\label{eq112}
    \frac{dz}{z}=\frac{dt}{t}.
\end{equation}
This shows that the corresponding similarity variable is
\begin{equation}\label{eq113}
    \eta=\frac{z}{t}.
\end{equation}
We now determine the transformation law for the dependent variable $v$. Along the characteristic curves, the third relation in \eqref{eq111} becomes
\begin{equation}\label{eq114}
    \frac{dv}{v}=-\left(\frac{\beta z^2}{4t^2}+\frac{1+\nu}{2t}\right)dt.
\end{equation}
Since $\eta=z/t$ remains constant along the characteristics, we may write $\eta^2=z^2/t^2$, and therefore integration gives
$$\ln v=-\frac{\beta}{4}\eta^2 t-\frac{1+\nu}{2}\ln t+\ln\varphi_5(\eta),$$
where $\varphi_5$ is an arbitrary function of the similarity variable $\eta$. Exponentiating both sides, we arrive at the invariant representation
\begin{equation}\label{eq115}
    v(z,t)=t^{-(1+\nu)/2}\exp\left(-\frac{\beta z^2}{4t}\right)\varphi_5\left(\frac{z}{t}\right).
\end{equation}
Thus, the invariant solution associated with $\widetilde{Y}_1$ is sought in the above similarity form. Substituting expression \eqref{eq115} into equation \eqref{eq108}, after straightforward calculation we obtain the reduced ordinary differential equation
\begin{equation}\label{eq116}
    \eta\varphi_5^{\prime\prime}+\nu\varphi_5^{\prime}=0.
\end{equation}
This equation can be integrated directly. Hence
\begin{equation}\label{eq117}
    \varphi_5(\eta)=\begin{cases}
        c_1+c_2\frac{\eta^{1-\nu}}{1-\nu},&\nu\neq 1,\\
        c_1+c_2\ln\eta,&\nu=1,
    \end{cases}
\end{equation}
where $c_1$ and $c_2$ are arbitrary constants. Consequently, we obtain the invariant solution
\begin{equation}\label{eq118}
    v(z,t)=\begin{cases}
        t^{-(1+\nu)/2}\exp\left(-\frac{\beta z^2}{4t}\right)\left(c_1+c_2\frac{(z/t)^{1-\nu}}{1-\nu}\right),&\nu\neq 1,\\
        t^{-1}\exp\left(-\frac{\beta z^2}{4t}\right)\left(c_1+c_2\ln\left(\frac{z}{t}\right)\right),&\nu=1.
    \end{cases}
\end{equation}
Finally, returning to the original dependent variable $u$ through the inverse transformation \eqref{eq109}, we obtain the corresponding invariant solution of the original nonlinear equation \eqref{eq2}.

\textbf{Generator $\widetilde{Y}_2$.} Analogously, if we consider $v=v(z,t)$ is invariant under the action of $\widetilde{Y}_2$, then the invariance condition
$\widetilde{Y}_2(v-v(z,t))\big|_{v=v(z,t)}=0.$
Then, it yields the first-order equation
\begin{equation}\label{eq119}
    \frac{z}{2}v_z+t v_t=0.
\end{equation}
This relation shows that the solution remains constant along the characteristic curves generated by $\widetilde{Y}_2$. Solving the corresponding characteristic system
$$\frac{dz}{z/2}=\frac{dt}{t},$$
we define similarity variable associated with the generator $\widetilde{Y}_2$ is $\eta=\frac{z}{t}$. It follows that the invariant solution can be sought in the form
\begin{equation}\label{eq120}
    v(z,t)=\varphi_6(\eta),\quad\quad \eta=\frac{z}{\sqrt{t}}.
\end{equation}
Substituting \eqref{eq120} into equation \eqref{eq108}, we reduce the original partial differential equation to an ordinary differential equation for the similarity profile $\varphi_6$ in the following form
\begin{equation}\label{eq121}
    \varphi_6^{\prime\prime}(\eta)+\left(\frac{\nu}{\eta}+\frac{\beta\eta}{2}\right)\varphi_6^{\prime}=0.
\end{equation}
Introducing the auxiliary function $w(\eta)=\varphi_6^{\prime}(\eta)$, one can estimate that its solution is readily obtained in the form 
\begin{equation}\label{eq122}
    w(\eta)=c_1\eta^{-\nu}\exp\left(-\frac{\beta\eta^2}{4}\right),
\end{equation}
where $c_1$ is an arbitrary constant. Hence, using substitution $w=\varphi_6^{\prime}$ again, we arrive at
\begin{equation}\label{eq123}
    \varphi_6(\eta)=c_2+c_1\int\eta^{-\nu}\exp\left(-\frac{\beta\eta^2}{4}\right)d\eta,
\end{equation}
where $c_2$ is another integration constant. Consequently, the invariant solution of the transformed equation corresponding to the generator $\widetilde{Y}_2$ is given by
\begin{equation}\label{eq124}
    v(z,t)=c_2+c_1\int \frac{z^{-\nu}}{t^{\frac{1-\nu}{2}}}\exp\left(-\frac{\beta z^2}{4t}\right)dz.
\end{equation}
Thus, the generator $\widetilde{Y}_2$ leads to a self-similar reduction with similarity variable $\eta=z/\sqrt{t}$. The resulting invariant solutions are described by a one-parameter integral family in the transformed variable $v$, and the corresponding solutions for the original variable u are recovered through the inverse mapping $u=J^{-1}(v)$. This type of reduction is typical for dilation symmetries and provides a natural radial similarity form for the heat equation in the constant case.

\textbf{Generator $\widetilde{Y}_3$.} This operator represents translation in time. Therefore, the invariant solutions associated with $\widetilde{Y}_3$ are precisely the stationary solutions, that is, solutions independent of the time variable.  

Indeed, if $v=v(z,t)$ is invariant under the action of $\widetilde{Y}_3$, then the invariance condition takes the simple form
\begin{equation}\label{eq125}
    v_t=0.
\end{equation}
Hence, $v$ depends only on the spatial variable $z$, and we may write
\begin{equation}\label{eq126}
    v(z,t)=\varphi_7(z).
\end{equation}
Substituting \eqref{eq126} into the linear radial heat equation \eqref{eq108}, we obtain the ordinary differential equation
\begin{equation}\label{eq127}
    \varphi_7^{\prime\prime}+\frac{\nu}{z}\varphi_7^{\prime}=0.
\end{equation}
Therefore, solution of the equation \eqref{eq127} takes the form
\begin{equation}\label{eq128}
    \varphi_7(z)=\begin{cases}
        C_0+C_1\frac{z^{1-\nu}}{1-\nu},& \nu\neq 1,\\
        C_0+C_1\ln z,&\nu=1,
    \end{cases}
\end{equation}
with arbitrary constants $C_0$ and $C_1$. Thus, the invariant solution of equation  associated with the generator $\widetilde{Y}_3$ is
\begin{equation}\label{eq129}
    u(z,t)=\begin{cases}
        J^{-1}\left( C_0+C_1\frac{z^{1-\nu}}{1-\nu}\right),&\nu\neq 1,\\
        J^{-1}\left(C_0+C_1\ln z\right),&\nu=1.
    \end{cases}
\end{equation}
Therefore, the generator $\widetilde{Y}_3$ yields the class of steady-state solutions. These invariant solutions do not evolve in time and describe equilibrium temperature distributions in the transformed variable $v$, as well as the corresponding stationary regimes for the original variable $u$.

\textbf{Generator $\widetilde{Y}_4$.} We next consider the infinitesimal generator $\widetilde{Y}_4=-\frac{J(u)}{K(u)}\partial u$. Recalling the transformation \eqref{eq107} and using $\partial_u=\frac{d v}{du}\partial_v=K(u)\partial_v$ because $\frac{dv}{du}=J^{\prime}(u)=K(u)$, we observe that the generator $\widetilde{Y}_4$ takes the simpler form
$$\widetilde{Y}_4=-v\partial_v,$$
which immediately reduces to 
$$-v=0.$$ 
Hence, the only invariant solution associated with this generator is $$v(z,t)=0.$$
Returning now to the original dependent variable $u$, we use the inverse transformation \eqref{eq109} and from $v=0$ we obtain
$$u(z,t)=J^{-1}(0).$$
which is a constant solution of the original nonlinear equation. It follows that the generator $\widetilde{Y}_4$ does not produce a nontrivial similarity reduction, but only yields the trivial constant invariant solution corresponding to the zero level of the transformed variable $v$. In other words, the action of $\widetilde{Y}_4$ scales the variable $v$ itself, and the only solution left unchanged by this scaling is $v=0$.

Consequently, in the original variable $u$, the invariant solution is the constant value determined implicitly by the relation $\int K(u)du=0$. If the function $J(u)=\int K(u)du$ is invertible, this constant is uniquely determined, otherwise, the invariant solution should be understood as any constant $u=u_0$ satisfying $J(u_0)=0$.

When $\nu=2$ the invariant solution $v=v(z,t)$ for infinitesimal generators $\widehat{Y}_1$, $\widehat{Y}_2$, $\widehat{Y}_3$ and $\widehat{Y}_6$ can be estimated analogously as in generators $\widetilde{Y}_1$, $\widetilde{Y}_2$, $\widetilde{Y}_3$ and $\widetilde{Y}_4$ and there are two additional generators are appeared.

\textbf{Generator $\widehat{Y}_4$ for $\nu=2$.}  
We now consider the infinitesimal generator $\widehat{Y}_4=t\partial_z-\left(\frac{\beta}{2}z+\frac{t}{z}\right)v\partial_v$. Let $v=v(z,t)$ be an invariant solution corresponding to $\widehat{Y}_4$. Then, by the invariance criterion, $v$ must satisfy $\widehat{Y}_4(v-v(z,t))\big|_{v=v(z,t)}=0$ which yields the first-order linear equation
\begin{equation*}
    tv_z+\left(\frac{\beta}{2}z+\frac{t}{z}\right)v=0.
\end{equation*}
Equivalently, this relation may be rewritten as
\begin{equation}\label{eq130}
    v_z=-\left(\frac{\beta z}{2t}+\frac{1}{z}\right)v.
\end{equation}
This is an ordinary differential equation with respect to the spatial variable $z$, while $t$ plays the role of a parameter. Integrating with respect to $z$, we obtain the invariant form of the solution. Indeed, solution of \eqref{eq130} is
\begin{equation}\label{eq131}
    v(z,t)=\frac{\varphi_8(t)}{z}\exp\left(-\frac{\beta z^2}{4t}\right),
\end{equation}
where $\varphi_8$ is an arbitrary function of $t$. 

We now substitute this expression into the governing equation \eqref{eq108} in order to determine the unknown function $\varphi_8$. Differentiating $v(z,t)$ with respect to $t$, differentiation with respect to $z$ and after one more differentiation with respect to $z$, we obtain
\begin{equation}\label{eq132}
    v_t=\frac{1}{z}\left(\varphi_8^{\prime}(t)+\frac{\beta z^2}{4t^2}\varphi_8(t)\right)\exp\left(-\frac{\beta z^2}{4t}\right),
\end{equation}
\begin{equation}\label{eq133}
    v_z=-\left(\frac{1}{z^2}+\frac{\beta}{2t}\right)\varphi_8(t)\exp\left(-\frac{\beta z^2}{4t}\right),
\end{equation}
\begin{equation}\label{eq134}
    v_{zz}=\left(\frac{\beta}{2tz}+\frac{\beta^2 z}{4t^2}+\frac{2}{z^3}\right)\varphi_8(t)\exp\left(-\frac{\beta z^2}{4t}\right).
\end{equation}
Thus, the equation \eqref{eq135} becomes
\begin{equation}\label{eq135}
    \varphi_8^{\prime}(t)+\frac{1}{2t}\varphi_8(t)=0.
\end{equation}
This is a first-order linear equation and can be integrated directly. Its general solution is
\begin{equation}\label{eq136}
    \varphi_8(t)=C_0 t^{-1/2}.
\end{equation}
where $C_0$ is an arbitrary constant. Substituting this expression back into \eqref{eq131} for $v(z,t)$, we finally obtain
\begin{equation}\label{eq137}
    v(z,t)=\frac{C_0}{z\sqrt{t}}\exp\left(-\frac{\beta z^2}{4t}\right).
\end{equation}
Returning to the original dependent variable $u$, we use the inverse transformation \eqref{eq109}. Thus, this symmetry reduction leads to an explicit exact solution expressed through the inverse of the integral transformation
$$u(z,t)=J^{-1}\left(\frac{C_0}{z\sqrt{t}}\exp\left(-\frac{\beta z^2}{4t}\right)\right).$$
The obtained solution has the structure of a Gaussian-type profile multiplied by the radial factor $z^{-1}$ and the temporal decay $t^{-1/2}$, which is consistent with the spherical geometry of the problem.

\textbf{Generator $\widehat{Y}_5$ for $\nu=2$.} We now consider the infinitesimal generator $\widehat{Y}_5=\partial_z-\frac{1}{z}v\partial_v$. Let $v=v(z,t)$ be an invariant solution associated with the generator $\widehat{Y}_6$. Then, by the invariance criterion, the function $v(z,t)$ must satisfy $\widehat{Y}_5(v-v(z,t))\big|_{v=v(z,t)}=0$. This condition leads to the first-order differential equation
\begin{equation}\label{eq138}
    v_z-\frac{1}{z}v=0.
\end{equation}
Depending on the sign convention adopted for the invariant branch, the corresponding invariant solution may be written in the form
\begin{equation}\label{eq139}
    v(z,t)=\frac{\varphi_9(t)}{z}.
\end{equation}
We now substitute expression \eqref{eq139} into the governing equation \eqref{eq108} in order to determine the unknown function $\varphi_9(t)$. Hence, we have
\begin{equation}\label{eq140}
    \beta\frac{\varphi_9^{\prime}(t)}{z}=\frac{2\varphi_9(t)}{z^3}+\frac{2}{z}\left(-\frac{\varphi_9(t)}{z^2}\right).
\end{equation}
The terms on the right-hand side cancel identically, so the reduced equation becomes
\begin{equation}\label{eq141}
    \beta\frac{\varphi_9^{\prime}(t)}{z}=0.
\end{equation}
Hence, $\varphi_9^{\prime}(t)=0$ which shows that $\varphi_9$ is constant. Therefore, we have $\varphi_9(t)=C_1$, where $C_1$ is an arbitrary constant. Substituting this result into \eqref{eq139}, we obtain the explicit invariant solution
\begin{equation}\label{eq142}
    v(z,t)=\frac{C_1}{z}.
\end{equation}
Finally, plugging to the inverse transformation \eqref{eq109} we can verify that invariant solution of the equation \eqref{eq2} takes the form
\begin{equation}\label{eq143}
    u(z,t)=J^{-1}\left(\frac{C_1}{z}\right),
\end{equation}
where $J(u)=\int K(u)du$. Thus, the generator $\widehat{Y}_5$ yields a simple but nontrivial invariant solution whose transformed variable $v$ depends only on the radial coordinate through the inverse power law $z^{-1}$.

\section{Particular cases of the coefficients $C(u)$ and $K(u)$}

\subsection{Power-type of the coefficients}
In this section, we examine a particularly important subclass of coefficient functions, namely the power-law constitutive relations
\begin{equation}\label{eq144}
    K(u)=k_0u^{m},\quad\quad C(u)=c_0u^{n},\quad\quad k_0>0,\quad c_0>0.
\end{equation}
where $m$ and $n$ are real parameters. Such expressions are widely used in nonlinear heat-transfer theory and diffusion-type models, since they provide a natural description of media whose thermal conductivity and heat capacity vary with the temperature according to algebraic laws. In addition to their physical relevance, power-law coefficients are especially convenient from the analytical point of view, because they allow the general symmetry classification obtained above to be written in explicit form.

For these choices of $K(u)$ and $C(u)$, several quantities appearing in the determining equations and in the finite transformations can be computed directly. In particular, the integral transformation

\begin{equation}\label{eq145}
    J(u)=\int K(u)du=\begin{cases}
        \frac{k_0}{m+1}u^{m+1},&m\neq 1,\\
        k_0\ln u,& m=-1,
    \end{cases}
\end{equation}
while the ratio of the coefficients becomes
\begin{equation}\label{eq146}
    \frac{C(u)}{K(u)}=\frac{c_0}{k_0}u^{n-m}.
\end{equation}
Hence, the symmetry structure depends directly on whether the exponents $m$ and $n$ coincide or not. If $n=m$, then $\frac{C(u)}{K(u)}$ is constant and the equation belongs to the linearizable class. If $n\neq m$, then $\frac{C(u)}{K(u)}$ is non-constant and the equation falls into the genuinely nonlinear class. Therefore, the power-law model naturally separates into two subcases, each corresponding to one of the main branches of the general symmetry classification.

\textbf{Constant-ratio case ($m=n$).} If we take $\frac{C(u)}{K(u)}=\frac{c_0}{k_0}=\beta$ constant and introducing $v=J(u)$ which satisfies the linear radial heat equation \eqref{eq108}, where $J$ defined by \eqref{eq145}. Then, taking inverse of the equation \eqref{eq107} for $u$, one can derive
\begin{equation}\label{eq147}
    u(z,t)=\begin{cases}
        \left(\frac{m+1}{k_0}v(z,t)\right)^{\frac{1}{m+1}},& m\neq 1,\\
        \exp\left(\frac{v(z,t)}{k_0}\right), & m=1.
    \end{cases}
\end{equation}
Based on the results obtained in Section 2.2, the invariant solution of the $u$ generated by infinitesimal generators are:
\begin{itemize}
    \item From the projective generator $\widetilde{Y}_1$ for $m\neq -1$, we have
    \begin{equation}\label{eq148}
        u(z,t)=\begin{cases}
            \left(\frac{m+1}{k_0}t^{-(1+\nu)/2}\exp\left(-\frac{\beta z^2}{4t}\right)\left(c_1+c_2\frac{(z/t)^{1-\nu}}{1-\nu}\right)\right)^{\frac{1}{m+1}}, & \nu\neq 1\\
            \left(\frac{m+1}{k_0 t}\exp\left(-\frac{\beta z^2}{4t}\right)\left(c_1+c_2\ln\left(\frac{z}{t}\right)\right)\right)^{\frac{1}{m+1}},&\nu=1
        \end{cases}
    \end{equation}
    and for $m=-1$, we get
    \begin{equation}\label{eq149}
        u(z,t)=\begin{cases}
            \exp\left(\frac{1}{k_0}t^{-(1+\nu)/2}\exp\left(-\frac{\beta z^2}{4t}\right)\left(c_1+c_2\frac{(z/t)^{1-\nu}}{1-\nu}\right)\right), &\nu\neq 1,\\
            \exp\left(\frac{1}{t}\exp\left(-\frac{\beta z^2}{4t}\right)\left(c_1+c_2\ln\left(\frac{z}{t}\right)\right)\right),&\nu=1.
        \end{cases}
    \end{equation}

    \item From the scaling generator $\widetilde{Y}_2$, the invariant solution of $u(z,t)$ can be defined as
    \begin{equation}\label{eq150}
        u(z,t)=\begin{cases}
            \left[\frac{m+1}{k_0}\left(c_2+c_1\int \frac{z^{-\nu}}{t^{\frac{1-\nu}{2}}}\exp\left(-\frac{\beta z^2}{4t}\right)dz\right)\right]^{\frac{1}{m+1}},& m\neq -1,\\
            \exp\left[\frac{1}{k_0}\left(c_2+c_1\int \frac{z^{-\nu}}{t^{\frac{1-\nu}{2}}}\exp\left(-\frac{\beta z^2}{4t}\right)dz\right)\right],& m=-1.
        \end{cases}
    \end{equation}

    \item From the time-translation reduction $\widetilde{Y}_3$, when $m\neq -1$, the invariant solution of the equation \eqref{eq2} becomes
    \begin{equation}\label{eq151}
        u(z,t)=\begin{cases}
            \left[\frac{m+1}{k_0}\left(C_0+C_1\frac{z^{1-\nu}}{1-\nu}\right)\right]^{\frac{1}{m+1}}, &\nu\neq 1,\\
            \left[\frac{m+1}{k_0}\left(C_0+C_1\ln z\right)\right]^{\frac{1}{m+1}},& \nu=1
        \end{cases}
    \end{equation}
    and for $m=-1$
    \begin{equation}\label{eq152}
        u(z,t)=\begin{cases}
            \exp\left[\frac{1}{k_0}\left(C_0+C_1\frac{z^{1-\nu}}{1-\nu}\right)\right],&\nu\neq 1,\\
            \exp\left[\frac{1}{k_0}\left(C_0+C_1\ln z\right)\right], &\nu=1.
        \end{cases}
    \end{equation}
    \item According to the last generator $\widetilde{Y}_4$, the invariant solution of the governing equation \eqref{eq2} is $u(z,t)=\text{const}$.
\end{itemize}

\textbf{Non constant ratio case ($n\neq m$).} In this case, the generated invariant solutions $u(z,t)$ by infinitesimal generators $Y_1$, $Y_2$, $Y_3$ and $Y_3$ are stated as following form:

\begin{itemize}
    \item The basic scaling invariant solution is defined by $u(\eta,t)=\varphi_1(\eta)$ where 
    \begin{equation}\label{eq153}
        \varphi_1(\eta)=C_2+C_1\int\frac{\eta^{-\nu}}{\varphi^{m}}\exp\left(-\frac{c_0}{2k_0}\int\eta\varphi^{n-m}(\eta)d\eta\right)d\eta, \quad \eta=\frac{z}{\sqrt{t}}.
    \end{equation}

    \item Since the generator $Y_2=\partial_t$ corresponds to invariance with respect to time translations, the associated invariant solutions are stationary. Therefore, the dependent variable does not depend on t, and we seek solutions in the form $u(z,t)=\varphi_2(z)$. In this case, $u_t=0$, for the power-law conductivity the equation \eqref{eq2} becomes
    \begin{equation}\label{eq154}
        \left(k_0z^{\nu}\varphi_2^{m}\varphi_2^{\prime}\right)^{\prime}=0.
    \end{equation}
    Integrating once with respect to $z$, we obtain
    $$k_0z^{\nu}\varphi_2^{m}\varphi_2^{\prime}=C_1,$$
    where $C_1$ is an integration constant. At this point, the explicit form of the stationary solution depends on the value of the exponent $m$. We therefore distinguish two cases.

    If $m\neq -1$, then integration of the left-hand side gives
    \begin{equation}\label{eq155}
        \frac{\varphi_2^{m+1}}{m+1}=\frac{C_1}{k_0}\int z^{-\nu}dz+C_2,
    \end{equation}
    where $C_2$ is another constant of integration. The form of the right-hand side depends on the parameter $\nu$. 
    
    If $m=-1$, then the equation \eqref{eq154} becomes
    \begin{equation}
        \ln\varphi_2=\frac{C_1}{k_0}\int z^{-\nu}dz+C_2.
    \end{equation}
    Therefore, for the special exponent $m=-1$, the invariant stationary solutions are expressed in logarithmic form for $\ln \varphi_2$, or, equivalently, in exponential form for $\varphi_2$ itself.
    
    Summarizing, the invariant solutions associated with the generator $Y_2$ are stationary and are given by
    \begin{equation}\label{eq156}
        u(z,t)=\begin{cases}
            \left(C_2+C_1 z^{1-\nu}\right)^{\frac{1}{m+1}},&m\neq -1,\quad \nu\neq 1,\\
            \left(C_2+C_1\ln z\right)^{\frac{1}{m+1}},& m=-1,\quad \nu=1
        \end{cases}
    \end{equation}
    and
    \begin{equation}\label{eq157}
        u(z,t)=\begin{cases}
            \exp\left(C_2+C_1z^{1-\nu}\right),&m=-1,\quad \nu\neq 1,\\
            C_2z^{C_1},&m=-1,\quad \nu=1.
        \end{cases}
    \end{equation}
    Thus, the symmetry generated by time translation leads to a family of exact stationary solutions whose form is determined explicitly by the power-law exponent $m$ and the geometric parameter $\nu$.

    \item The existence of the generator $Y_3$ is governed by the function \eqref{eq20} and for the present power-law choice, it takes the form $F(u)=\frac{u}{n-m}$. On the other hand, according to the general result in Section 2, the existence of the additional generator $Y_3$ requires that $F(u)$ can be represented in the form \eqref{eq38}. Hence, to determine whether the extra symmetry is present for power-law coefficients, it is enough to equalize \eqref{eq38} with the explicit power-law formula for $F(u)$ is defined by \eqref{eq20}. For $m=-1$, we have
    \begin{equation}\label{eq159}
        \frac{k_0}{n-m}u^{m+1}=\frac{Ak_0}{m+1}u^{m+1}+B.
    \end{equation}
    Since this identity must hold for all admissible values of u, the constant term and the coefficient of $u^{m+1}$ must coincide separately. Therefore, we must have
    $$A=\frac{m+1}{n-m},\quad\quad B=0.$$
    Hence, for $n\neq m$ and $m\neq -1$, the required compatibility condition is satisfied automatically, and therefore the additional symmetry $Y_3$ does indeed exist in the power-law case.

    We now specialize the general invariant formula associated with the generator $Y_3$. From the previous Section 2, the invariant representation corresponding to this symmetry can be written as
    $$AJ(u)=z^{-2A}\left[Q+\frac{2(2A+1-\nu)}{D}t\right]^{A},$$
    where $Q$ and $D$ are arbitrary constants. Substituting $J(u)=\frac{k_0}{m+1}u^{m+1}$ and $A=\frac{m+1}{n-m}$, we obtain
    $$\frac{k_0}{n-m}u^{m+1}=z^{-2A}\left[Q+\frac{2(2A+1-\nu)}{D}t\right]^{A}.$$
    Solving for $u(z,t)$, we find
    \begin{equation}\label{eq160}
        u(z,t)=\left[\frac{n-m}{k_0}z^{-2A}\left(Q+\frac{2(2A+1-\nu)}{D}t\right)^{A}\right]^{\frac{1}{m+1}},\quad\quad A=\frac{m+1}{n-m}.
    \end{equation}
    Therefore, the invariant solution associated with the extra generator $Y_3$ has the form of a separable power law in the radial variable $z$, multiplied by a time-dependent factor determined by the symmetry parameters and the constants of integration. This result shows that the additional scaling symmetry admitted by the power-law coefficients leads to an explicit family of exact similarity solutions.
\end{itemize}

\subsection{Exponential-type of the coefficients}
We now consider another physically relevant subclass of constitutive functions, namely the exponential-type laws
\begin{equation}\label{eq161}
    K(u)=k_0e^{\lambda u},\quad\quad C(u)=c_0e^{\mu u},\quad\quad k_0>0,\quad c_0>0,
\end{equation}
where $\lambda,\mu\in\mathbb{R}$ are constants. Such coefficients arise naturally in nonlinear heat-transfer models when the thermal conductivity and heat capacity vary exponentially with the temperature. From the point of view of symmetry analysis, this class is also especially convenient, since the integral transformation $J(u)=\int K(u)du$ and the quotient $\frac{C(u)}{K(u)}$ can be computed explicitly.

Indeed,
\begin{equation}\label{eq162}
    J(u)=\int K(u)du=\begin{cases}
        \frac{k_0}{\lambda}e^{\lambda u},&\lambda\neq 0,\\
        k_0u,&\lambda=0
    \end{cases}
\end{equation}
and
\begin{equation}\label{eq163}
    \frac{C(u)}{K(u)}=\frac{c_0}{k_0}e^{(\mu-\lambda)u}.
\end{equation}
Therefore, just as in the power-law case, the symmetry structure is determined by whether the ratio $C(u)/K(u)$ is constant or not. This leads naturally to two subcases, if $\mu=\lambda$, then $C(u)/K(u)=c_0/k_0=\beta$ constant, and the equation belongs to the constant-ratio class. In this case, the substitution $v=J(u)$ transforms the original nonlinear equation \eqref{eq2} into the linear radial heat equation \eqref{eq108}. On the other hand, if $\mu\neq\lambda$, then the ratio $C(u)/K(u)$ remains genuinely dependent on $u$, and the equation falls into the non-constant-ratio class. We treat these two cases separately below.

\textbf{Case $\mu=\lambda$.} Thus, all invariant solutions obtained earlier for the constant-ratio case can be transferred directly to the original variable $u$ through the inverse transformation

\begin{equation}\label{eq164}
    u(z,t)=\begin{cases}
        \frac{1}{\lambda}\ln\left(\frac{\lambda}{k_0}v(z,t)\right),&\lambda\neq 0,\\
        \frac{v(z,t)}{k_0},&\lambda=0.
    \end{cases}
\end{equation}
\begin{itemize}
    \item $\widetilde{Y}_1$:  For $\lambda\neq 0$, we have the following invariant solution form
    \begin{equation}\label{eq165}
        u(z,t)=\begin{cases}
            \frac{1}{\lambda}\ln\left[\frac{\lambda}{k_0}t^{-(1+\nu)/2}\exp\left(-\frac{\beta z^2}{4t}\right)\left(c_1+c_2\frac{(z/t)^{1-\nu}}{1-\nu}\right)\right],&\nu\neq 1,\\
            \frac{1}{\lambda}\ln\left[\frac{\lambda}{k_0}t^{-1}\exp\left(-\frac{\beta z^2}{4t}\right)\left(c_1+c_2\ln\left(\frac{z}{t}\right)\right)\right],&\nu=1.
        \end{cases}
    \end{equation}
    Analogously, we can obtain for $\lambda=0$.

    \item $\widetilde{Y}_2$: According to the generator $\widetilde{Y}_2$, the invariant solution takes the form
    \begin{equation}\label{eq166}
        u(z,t)=\begin{cases}
            \frac{1}{\lambda}\ln\left[\frac{\lambda}{k_0}\left(c_2+c_1\int \frac{z^{-\nu}}{t^{\frac{1-\nu}{2}}}\exp\left(-\frac{\beta z^2}{4t}\right)dz\right)\right],&\lambda\neq 0,\\
            \frac{1}{k_0}\left[c_2+c_1\int \frac{z^{-\nu}}{t^{\frac{1-\nu}{2}}}\exp\left(-\frac{\beta z^2}{4t}\right)dz\right],&\lambda=0.
        \end{cases}
    \end{equation}

    \item $\widetilde{Y}_3$: The invariant solution corresponding to time translation is obtained by setting $v_t=0$. Based on \eqref{eq129} we have the following solutions
    \begin{equation}\label{eq167}
        u(z,t)=\begin{cases}
            \frac{1}{\lambda}\ln\left[\frac{\lambda}{k_0}\left( C_0+C_1\frac{z^{1-\nu}}{1-\nu}\right)\right],&\lambda\neq 0,\quad\nu\neq 1,\\
            \frac{1}{\lambda}\ln\left[\frac{\lambda}{k_0}\left(C_0+C_1\ln z\right)\right],&\lambda=0,\quad\nu=1.
        \end{cases}
    \end{equation}
    Applying \eqref{eq128} and \eqref{eq164} one can get invariation solution for $\lambda=0$, respectively.
\end{itemize}

\textbf{Case $\mu\neq\lambda$.} We now turn to the genuinely nonlinear case, in which invariant function $u$ of the original nonlinear generalized heat-diffusion equation \eqref{eq2}, admitted by the infinitesimal generators \eqref{eq47} can be determined as follows.

\begin{itemize}
    \item $Y_1$: The similarity solution can be verified as $u(z,t)=\varphi_1(\eta)$ with $\eta=z/\sqrt{t}$ and $\varphi_1$ takes the form
    \begin{equation}\label{eq168}
        \varphi_1(\eta)=C_2+\frac{C_1}{k_0}\int \eta^{-\nu}\exp\left(-\lambda\varphi_1-\frac{c_0}{2k_0}\int \eta e^{(\mu-\lambda)\varphi_1}d\eta\right)d\eta.
    \end{equation}

    \item $Y_2$: By this generator, the invariant solution is presented as $u(z,t)=\varphi_2(z)$, where
    \begin{equation}\label{eq169}
        \varphi_2(z)=\begin{cases}
            \frac{1}{\lambda}\ln\left(\frac{\lambda C_1}{k_0(1-\nu)}z^{1-\nu}+C_2\right),&\nu\neq 1,\\
            \frac{1}{\lambda}\ln\left(\frac{\lambda C_1}{k_0}\ln z+C_2\right), &\nu=1.
        \end{cases}
    \end{equation}

    \item $Y_3$: Substituting the function $F(u)$ defined by \eqref{eq20} into the equation \eqref{eq38} we have
    $$\frac{1}{\mu-\lambda}=\frac{A}{\lambda}+\frac{B}{k_0}e^{-\lambda u},$$
    which leads to $A=\frac{\lambda}{\mu-\lambda}$ and $B=0$. Then invariant relation \eqref{eq89} is represented in the form 
    $$\frac{k_0}{\mu-\lambda}e^{\lambda u}=z^{-2A}\left(Q+\frac{2(2A+1-\nu)}{D}t\right)^A,$$
    when $\lambda\neq 0$. Hence, solving for $u(z,t)$, we have
    \begin{equation}\label{eq170}
        u(z,t)=\frac{1}{\lambda}\ln\left[\frac{\mu-\lambda}{k_0}z^{-2A}\left(Q+\frac{2(2A+1-\nu)}{D}t\right)^A\right],
    \end{equation}
    with $D$ and $Q$ are arbitrary constants.
\end{itemize}

\subsection{Linear case of the coefficients}
We now consider the case when the conductivity and heat capacity depend linearly on the temperature:
\begin{equation}\label{eq171}
    K(u)=k_0(a+bu),\quad\quad C(u)=c_0(c+du),\quad\quad k_0>0,\quad c_0>0,
\end{equation}
where $a,b,c,d\in\mathbb{R}$. This class is of interest because it represents the simplest nontrivial departure from constant coefficients, while still allowing the principal quantities arising in the Lie symmetry classification to be written explicitly. For these coefficients, the integral transformation
\begin{equation}\label{eq172}
    J(u)=\int K(u)du=\begin{cases}
        k_0\left(au+\frac{bu^2}{2}\right),& b\neq 0,\\
        k_0au,&b=0.
    \end{cases}
\end{equation}
Moreover, the key ratio entering the determining equations is
\begin{equation}\label{eq173}
    \frac{C(u)}{K(u)}=\frac{c_0(c+du)}{k_0(a+bu)}.
\end{equation}
Hence, invariant solution function $u(z,t)$ for governing equation \eqref{eq2} can be verified:

\begin{itemize}
    \item with infinitesimal generator $Y_1$ in the form of $u(z,t)=\varphi_1(\eta)$ where $\eta=z/\sqrt{t}$ and
    \begin{equation}\label{eq174}
        \varphi_1(\eta)=C_2+\frac{C_1}{k_0}\int\frac{\eta^{-\nu}}{a+b\varphi_1}\exp\left(-\frac{c_0}{2k_0}\int\eta\frac{c+d\varphi_1}{a+b\varphi_1}d\eta\right)d\eta.
    \end{equation}

    \item with generator $Y_2$ and takes the form $u(z,t)=\varphi_2(z)$ where $\varphi_2$ can be determined implicitly from the equations
    \begin{equation}\label{eq175}
        a\varphi_2+\frac{b\varphi_2^2}{2}=\frac{C_1}{k_0(1-\nu)}z^{1-\nu}+C_2,\quad\quad \nu\neq 1
    \end{equation}
    and
    \begin{equation}\label{eq176}
        a\varphi_2+\frac{b\varphi_2^2}{2}=\frac{C_1}{k_0}\ln z+C_2,\quad\quad \nu=1.
    \end{equation}

    \item by generator $Y_3$, from the implicit equation
    \begin{equation}\label{eq177}
        au+\frac{bu^2}{2}=\frac{1}{k_0 A(u)}z^{-2A(u)}\left(Q+\frac{2(2A(u)+1-\nu)}{D}t\right)^{A(u)},
    \end{equation}
    where $A(u)=\frac{d(a+bu)}{(au+\frac{bu^2}{2})(c+du)}-\frac{b}{au+\frac{bu^2}{2}}$.
    \item with generator $Y_4$ which implies $u(z,t)\equiv u_0$.
\end{itemize}

\section{Conclusion}
In this paper, we investigated the nonlinear generalized heat equation which models heat conduction in a medium with temperature-dependent heat capacity and thermal conductivity in a variable cross-section geometry determined by the parameter $\nu>0$. The main objective was to classify the admitted Lie point symmetries of this equation and to use them for symmetry reduction and construction of invariant solutions.

The symmetry analysis showed that the classification is governed by the functional structure of the ratio $C(u)/K(u)$. This leads naturally to two principal branches: the genuinely nonlinear case when $C(u)/K(u)$ is non-constant, and the proportional case $C(u)=\beta K(u)$, where the equation becomes reducible through the transformation $v=J(u)=\int K(u)du$ to the linear radial heat equation \eqref{eq108}. For the non-constant-ratio case, the equation always admits the basic generators corresponding to joint space-time scaling and time translation, while additional generators appear only under special compatibility relations between $C(u)$ and $K(u)$. For the constant-ratio case, a richer symmetry structure was obtained, with a four-dimensional Lie algebra for $\nu\neq 2$ and a six-dimensional Lie algebra in the spherical case $\nu=2$.

On the basis of these admitted generators, a family of invariant and similarity solutions was constructed. In the non-constant branch, the reductions associated with the generators $Y_1$, $Y_2$, $Y_3$ and $Y_4$ led to reduced ordinary differential equations and, in several cases, to explicit invariant representations. In the constant-ratio branch, the transformation to the linear radial heat equation made it possible to derive invariant solutions systematically in the transformed variable $v$, and then recover the corresponding solutions for the original variable $u$ through the inverse mapping $J^{-1}$. In the spherical case $\nu=2$ was shown to admit additional symmetries $\widehat{Y}_4$, $\widehat{Y}_5$ and $\widehat{Y}_6$.

The general classification was then specialized to several important subclasses of constitutive laws. For power-law coefficients, the analysis showed explicitly how the symmetry structure depends on whether the exponents coincide or not, and explicit invariant solutions were obtained in both the constant and non-constant ratio regimes. Exponential-type conductivity/capacity pairs were also examined, again illustrating both branches of the general classification and providing closed-form invariant solutions in representative cases. These examples demonstrate that the abstract symmetry framework developed in the paper can be applied effectively to concrete nonlinear heat-transfer models with physically meaningful constitutive behavior.

Overall, the results show that Lie symmetry analysis provides a unified and effective tool for studying generalized nonlinear heat equations with variable geometry. The obtained infinitesimal generators, commutator structures, one-parameter groups, and invariant solutions form a basis for further analytical study of nonlinear diffusion-type models in radial, cylindrical, and spherical geometries.

\textbf{Limitations and future perspectives.}

The limitation of the study is that the complete symmetry enhancement in the non-constant-ratio case occurs only under special compatibility relations between the constitutive functions $C(u)$ and $K(u)$. Accordingly, for arbitrary nonlinear coefficients the admitted algebra is relatively small, and the resulting invariant reductions may lead only to implicit integral relations rather than fully explicit solutions. Thus, the current classification is strongest for special functional subclasses and less explicit for completely general constitutive laws.

Another limitation is that the paper focuses on the differential equation itself and does not include boundary and initial conditions. From the physical point of view, however, many heat-conduction problems of interest involve free boundaries, moving interfaces, or prescribed boundary fluxes. Some of the admitted symmetries obtained here may not be compatible with such constraints, so the invariant solutions derived in the paper should be interpreted primarily as exact solutions of the governing equation rather than as complete solutions of specific boundary-value problems. This is especially relevant for applications to Stefan-type and phase-change problems.

The one of the important perspective is to apply the present classification to boundary-value and free-boundary problems. Since the generalized heat equation studied here is closely connected with radial diffusion and phase-change models, it would be valuable to identify which of the admitted symmetries remain compatible with physically meaningful boundary conditions and to use them for constructing invariant solutions of Stefan-type problems in cylindrical and spherical geometries. This would connect the abstract group classification more directly with applications in thermal analysis and electrical-contact modeling.

The one of the important perspective is to apply the present classification to boundary-value and free-boundary problems. Since the generalized heat equation studied here is closely connected with radial diffusion and phase-change models, it would be valuable to identify which of the admitted symmetries remain compatible with physically meaningful boundary conditions and to use them for constructing invariant solutions of Stefan-type problems in cylindrical and spherical geometries. This would connect the abstract group classification more directly with applications in thermal analysis and electrical-contact modeling.

A further direction is to combine the Lie symmetry reductions obtained here with integral methods, similarity techniques, and inverse-problem approaches in order to study more complicated models involving source terms, variable geometry, or coupled thermal effects. Such developments could significantly extend the applicability of the present results to realistic models in heat transfer, phase transitions, and contact phenomena.

\end{document}